\documentclass{article}
\usepackage{graphicx}
\usepackage{amsmath, amssymb}
\usepackage{amsthm}
\usepackage{amsfonts}
\usepackage{bm}

\usepackage{tikz}
\usepackage{caption}

\usepackage{tikz, pgfplots}
\pgfplotsset{compat=1.18}

\usepackage{titlesec}
\titleformat{\section}{\large\bfseries}{\thesection.}{1em}{}

\begin{document}

Tsemo Aristide

PKFokam Institute of Technology

Yaounde, Cameroon.

tsemo58@yahoo.ca

\bigskip
\bigskip

\centerline{\bf The algebra and the geometry aspect of Deep learning.}

\bigskip
\bigskip

\begin{abstract}
This paper investigates the foundations of deep learning through insight of geometry, algebra and differential calculus. At is core, artificial intelligence relies on assumption that data and its intrinsic structure can be embedded into vector spaces allowing for analysis through geometric and algebraic methods. We thrace the development of neural networks from the perceptron to the transformer architecture, emphasizing on the underlying geometric structures and differential processes that govern their behavior. Our original approach highlights how the canonical scalar product on matrix spaces naturally leads to backpropagation equations yielding to a coordinate free formulation.  We explore how classification problems can reinterpreted using tools from differential and algebraic geometry suggesting that manifold structure, degree of variety, homology may inform both convergence and interpretability of learning algorithms We further examine how neural networks can be interpreted via their associated directed graph, drawing connection to a Quillen model defined in [1] and [13] to describe memory as an homotopy theoretic property of the associated network.
\end{abstract}

Most of the algorithms refered in this paper are available at:

 https://github.com/TsemoAristide/Deep-learning-explained-to-a-mathematician

\section*{1. Introduction.}

Machine learning is a branch of artificial intelligence dedicated to developing algorithms capable of performing tasks after being properly trained. These algorithms are generally classified into three main categories:

\begin{itemize}
    \item \textbf{Supervised learning}, where the algorithm learns patterns from labeled training data;
    \item \textbf{Unsupervised learning}, where it identifies structures or groupings in unlabeled data;
    \item \textbf{Reinforcement learning}, where the algorithm improves its performance through trial-and-error interactions with an environment.
\end{itemize}

Within supervised learning, one important category is {classification algorithms}. These are designed to determine the class or category of a given data point, based on a partition of a subset \( S \) of a vector space. Once trained, the algorithm can classify new elements from the same space.
Such algorithms are well-suited to data that can be faithfully represented as vectors  such as images, sounds, or texts.

One of the earliest classification algorithms is the {perceptron}, introduced in [11]. It is a \textit{linear classifier}, meaning it can identify a separating hyperplane between two finite subsets of a vector space  if such a hyperplane exists.
A more sophisticated version of the perceptron is {logistic regression}. Like the perceptron, it can separate two classes, but it also estimates the {probability} that a given element belongs to a specific class. The training process involves {minimizing a loss function} typically the  cross-entropy.

However, some subsets of vector spaces are {separable}, but not \textit{linearly separable}. In such cases, more advanced classification algorithms  such as support vector machines with kernels or neural networks  are required.

When such datasets are embedded in low-dimensional vector spaces, it becomes possible to visualize them by plotting, allowing us to observe the shape of the separating surfaces and to design suitable classification algorithms accordingly. This approach is echoed in traditional classification methods, such as the Viola--Jones algorithm [15] in computer vision, where predefined patterns (features) are combined with machine learning techniques to perform classification.

{Deep learning} is a modern branch of machine learning capable of handling the classification of such datasets without the need for manually defined features. In this context, the output function that determines class probabilities is defined by a {neural network}.

Although machine learning can be highly efficient, allowing it to replace humans in the decision-making process can be risky.

Algorithms are trained on data collected up to a specific point in time. If the nature of a disease changes or new variants emerge (like we saw with COVID-19), the algorithm may no longer be effective. This phenomenon is called data drift, where the distribution of the data the model was trained on no longer reflects the current reality. Without continuous retraining or adaptation, the model's performance can degrade, leading to misclassifications or missed diagnoses.

The  decision-making processes of many deep learning models are not easily interpretable by humans. This lack of transparency can be dangerous in high-stakes fields like healthcare, where understanding why a model makes a certain prediction is crucial for trust and accountability.

The {universal approximation theorem} states that neural networks can approximate any continuous function on a compact subset of a vector space to an arbitrary degree of accuracy. While other approximation techniques, such as {Fourier series} can  achieve similar goals, neural networks have the unique advantage of {automatically discovering classifying patterns}, which can sometimes be interpreted. This capacity is particularly significant in the context of \textbf{generative artificial intelligence}, where understanding and synthesizing complex patterns is crucial and has applications to art with algorithms like style transfer (see [6]) where the network learns and reimagines visual features from one image in the stylistic framework of another.

The training set typically represents a specific class of objects, such as human faces. It is important that a neural network trained to classify whether an image contains a human face can generalize successfully to new, unseen images. From a mathematical perspective, it is always possible to find a separating manifold for a finite dataset. However, such a solution may not generalize well to unseen data, a problem known as {overfitting}.

To mitigate overfitting, it is often desirable to find separating manifolds that correctly classify a high percentage of the training dataset, but not necessarily all elements. According to \textbf{Nash's Theorem}, every compact manifold embedded in a Euclidean space can be approximated arbitrarily closely by an algebraic variety. Thus, we may assume that the classifying manifold is an {algebraic variety}. Given a target classification accuracy, the optimal separating manifold can be considered the algebraic variety of {lowest degree} that achieves this level of accuracy.

\bigskip

\noindent\textbf{More precisely:}

Let \( S \) be a finite subset of a vector space, and let
\[
f: S \rightarrow \{0, 1\}
\]
be a binary labelling function. Define the subsets:
\[
S_0 = f^{-1}(0), \quad S_1 = f^{-1}(1).
\]
Let \( p \in [0, 1] \) represent the desired classification accuracy. The goal is to find an {algebraic hypersurface} (i.e., the zero set of a polynomial) of {minimal degree} that separates \( S_0 \) and \( S_1 \), such that at least a proportion \( p \) of the elements in \( S \) are correctly classified. This is an aspect of topological data analysis (see [5]).

In deep learning, the more data one has, the more accurate the results tend to be. However, training on large datasets is often very costly in terms of both time and computational resources. From this perspective, finding a separating algebraic manifold of the lowest possible degree may be beneficial, as it could offer a more efficient way to achieve acceptable classification accuracy with less data or simpler models. This leads to the following conjecture:

{\bf CPU-GPU conjecture.}
Any result achieved by a model trained with  graphical processing units (GPU) can be achieved equivalently by a model trained with a central processing unit (CPU) in comparable time, provided that the CPU trained model incorporates sufficient mathematical insight.

Deep learning models such as EfficientNet, ResNet,... have been widely used to address a variety of tasks. These models are typically optimized through trial-and-error methods to identify the best hyperparameters, a process that can be both time-consuming and computationally expensive.

A common strategy for working with new datasets is fine-tuning, where a model pre-trained on a related dataset is adapted to the new task by modifying its final layers. While this approach is widely used and often effective, the nature of the decision boundaries (or separating manifolds) formed by fine-tuned models remains poorly understood. Gaining deeper insight into these structures could enable more efficient to design new model  for related problems, reducing the reliance on costly trial-and-error methods.

On the other hand, artificial intelligence has already had a significant impact on fundamental scientific research, for example, in solving the protein folding problem in 2020 by {DeepMind}. However, in mathematics, AI has not yet succeeded in resolving fundamental conjectures. (See [7]).

We propose the development of a \textbf{mathematical language model} aimed at studying and solving mathematical problems. The training set would consist of all known propositions and theorems, tokenized and learned using next-token prediction, as is done in natural language processing.
With appropriate alignment algorithms, such a model could potentially solve new mathematical problems provided their solutions rely on methods similar to those used in previously solved problems included in the training data.
The main obstacle to building such a model lies in the construction of its dataset. Unlike general language models like GPT, which are trained on the vast corpus of data available on the internet, mathematical knowledge particularly in formal or semi-formal form is far less accessible and standardized.
Remark that a similar structure exists for coding with the model CodeGen which completes codes.

The aim of this note is to present the mathematical foundations of deep learning and to implement key concepts from scratch using Python. 

We start by reviewing the fundamental principles of optimization, which are central to this field. Since datasets are often represented as matrices, we also cover differential calculus on matrices. We use the standard scalar product defined on a space of matrices defined by $\langle A,B\rangle =tr(B^TA)$ and Jacobian matices. This original approach simplifies the traditional presentation of gradient descent in deep learning.

We continue by introducing the perceptron algorithm, followed by a discussion on logistic regression. We then explore artificial neural networks (ANNs), which are capable of solving non-linear classification tasks. This leads us to the {Universal Approximation Theorem}, which states that any classification task can, in theory, be accomplished using ANNs. For theoretical rigorousness we view the rectified linear unit as a distribution to compute its derivative.

However, while ANNs are theoretically powerful, they can be computationally expensive, especially for large datasets. To address this, we study specialized classes of neural networks. We describe {Convolutional Neural Networks (CNNs)}, which are well-suited for image data. We give two applications of CNN to health. We construct with Pytorch a classifying network to detect malaria. With Tensorflow, we use the technique of transfer learning to build an algorithm which detects pneumonia among children. 

After, we study {Recurrent Neural Networks (RNNs)}, which are designed for sequential data.
We conclude this section with a presentation of the {Transformer architecture}, the core concept behind large language models such as {ChatGPT}, {Ollama}, and {DeeSeek}.

We also study the underlying architecture of  neural networks. They are objects of the category of directed graphs. We relate the class of the architecture of a neural network in the closed model defined in [13], to its memory. We study graph neural networks and endow the category of directed graphs with a tensor product.

Most of the algorithms discussed have been implemented in Python and are available in the accompanying repository available at:

 https://github.com/TsemoAristide/Deep-learning-explained-to-a-mathematician.

\bigskip

{\bf Plan.}

\medskip

0. Introduction.

1. Preliminaries.

2. The Perceptron.

3. Logistic Regression

4. Deep Neural Networks.

5. Universal Approximation Theorem.

6. Overfitting.

9. The category of Neural Networks.

7. Convolutional Neural Networks.

8. Reccurent Neural Networks and transformer

\section*{1. Preliminaries}

\subsection*{1.1. Optimization.}

\bigskip

Let \( \mathbb{R}^n \) be the \( n \)-dimensional real vector space endowed with its canonical basis.

Let \( u, v \in \mathbb{R}^n \) have coordinates \( (x_1, \dots, x_n) \) and \( (y_1, \dots, y_n) \), respectively. The bilinear form defined on \( \mathbb{R}^n \) by

\[
\langle u, v \rangle = x_1 y_1 + \dots + x_n y_n
\]

is positive definite. This is the \emph{canonical scalar product} associated with the standard basis of \( \mathbb{R}^n \).

For any element \( x \in \mathbb{R}^n \), the norm of \( x \) is defined as:

\[
\|x\| = \sqrt{\langle x, x \rangle}
\]

\bigskip

A function \( f: \mathbb{R}^n \rightarrow \mathbb{R} \) is said to be differentiable at a point \( x \in \mathbb{R}^n \) if there exists a linear map \( df_x: \mathbb{R}^n \rightarrow \mathbb{R} \) such that:

\[
f(x + h) = f(x) + df_x(h) + o(\|h\|)
\]

where:

\[
\lim_{\|h\| \to 0} \frac{o(\|h\|)}{\|h\|} = 0
\]

\medskip

{\bf 1.1.1 Gradient Descent.}

\medskip

Suppose that \( f \) is differentiable at the point \( x \in \mathbb{R}^n \). The gradient \( \nabla f(x) \) of \( f \), with respect to the scalar product \( \langle \cdot, \cdot \rangle \), is the vector in \( \mathbb{R}^n \) defined by:

\[
\langle \nabla f(x), u \rangle = df_x(u), \quad \text{for all } u \in \mathbb{R}^n.
\]

Suppose now that \( x \) is not a critical (extremum) point, i.e., \( \nabla f(x) \neq 0 \). Let \( t \in \mathbb{R} \). Then:

\begin{align*}
f(x + t \nabla f(x)) 
&= f(x) + df_x(t \nabla f(x)) + o(\|t \nabla f(x)\|) \\
&= f(x) + t \|\nabla f(x)\|^2 + o(t \|\nabla f(x)\|).
\end{align*}

This implies:

\[
f(x + t \nabla f(x)) - f(x) = t \left( \|\nabla f(x)\|^2 + \|\nabla f(x)\| \cdot \frac{o(t \nabla f(x))}{t \|\nabla f(x)\|} \right).
\]

\medskip

Since
\[
\lim_{t \to 0} \frac{o(t \nabla f(x))}{t \|\nabla f(x)\|} = \lim_{t \to 0} \frac{o(t \nabla f(x))}{\|t \nabla f(x)\|} = 0,
\]
we deduce that there exists \( \varepsilon > 0 \) such that for all \( 0 < t < \varepsilon \), we have:

\[
f(x - t \nabla f(x)) < f(x).
\]

This shows that moving in the direction opposite to the gradient decreases the value of the function  this is the principle of gradient descent. This method is effective when it is not 
possible to solve the equation  \(df (x)=0\).

\medskip

\noindent\textbf{Gradient Descent Algorithm to Find a Local Minimum:}

\begin{itemize}
  \item Choose an initial point \( x_0 \in \mathbb{R}^n \).
  \item Choose a learning rate \( \alpha > 0 \).
  \item Repeat for \( t = 0, 1, \dots, N-1 \):
  \[
  x_{t+1} = x_t - \alpha \nabla f(x_t).
  \]
\end{itemize}

If \(\alpha\) is small enough, the sequence \(x_t\) will converge towards a local minimum.

This algorithm is implemented in the accompanying repository.

\medskip

{\bf 1.1.2. Momentum.}

\medskip

Momentum is an optimization technique used to accelerate gradient descent by incorporating information from past updates. Instead of relying solely on the current gradient, momentum maintains a moving average of past gradients to build up speed in relevant directions and dampen oscillations.

The update rules are defined as follows:

\begin{align*}
v_{t+1} &= \gamma v_t + \alpha \nabla f(x_t), \\
x_{t+1} &= x_t - v_{t+1},
\end{align*}

where:
\begin{itemize}
    \item $v_t$ is the \textbf{velocity} at time step $t$ and $v_0=0$,
    \item $\gamma \in [0, 1)$ is the \textbf{momentum coefficient} (typically $\gamma = 0.9$ or $0.99$),
    \item $\alpha$ is the \textbf{learning rate},
    \item $\nabla f(x_t)$ is the gradient of the objective function at $x_t$.
\end{itemize}

The key idea is that $v_t$ acts like a memory of past gradients. The term $\gamma v_t$ helps retain the previous direction of motion, while $\alpha \nabla f(x_t)$ incorporates the new gradient information. This combination can help the algorithm maintain consistent movement along shallow but long valleys in the loss surface and reduce erratic updates in directions where the gradient changes rapidly.

\medskip

\noindent \textbf{Advantages of Momentum:}
\begin{itemize}
    \item Speeds up convergence, especially in scenarios with high curvature or noisy gradients.
    \item Helps to escape shallow local minima or flat regions of the loss function.
    \item Reduces oscillations in the gradient descent path, particularly in directions where the gradient sign alternates.
\end{itemize}

\medskip

Overall, momentum introduces inertia into the optimization process, smoothing the trajectory of parameter updates and often leading to faster and more stable convergence compared to vanilla gradient descent.

\medskip

{\bf 1.1.3. RMSProp.}

\medskip

An important challenge in optimization is choosing an appropriate learning rate. RMSProp (Root Mean Square Propagation) is an adaptive learning rate algorithm that addresses this issue by adjusting the learning rate at each time step.

Let \( g_t= \nabla f(x_t) \) denote the gradient of the function \( f \) at step \( t \). RMSProp maintains an exponential moving average of the squared gradients:

\[
E_0 = 0
\]
\[
E_{t+1} = \beta E_t + (1 - \beta) g_t^2
\]

Here, \( g_t^2 \) denotes the element-wise square of the gradient vector \( g_t \), and \( \beta \in [0,1) \) is a decay rate hyperparameter, typically around 0.9.

The parameter update rule is given by:

\[
x_{t+1} = x_t - \frac{\eta}{\sqrt{E_{t+1} + \epsilon}} \cdot g_t
\]

where:
\
  \( \eta \) is the learning rate,
  \( \epsilon \) is a small constant (e.g., \( 10^{-8} \)) added for numerical stability,
  The division is element-wise.

This adaptive scheme ensures that:

   When the gradient components are large, the corresponding learning rates shrink, helping to prevent exploding gradients.
   When gradient components are small, the effective learning rates remain larger.

\medskip

{\bf 1.1.4. ADAM: Adaptive Moment Estimation.}

\medskip

ADAM is an algorithm that combines ideas from momentum and RMSProp.

\medskip

\textbf{Momentum term:}
\[
m_{t+1} = \beta_1 m_t + (1 - \beta_1) \nabla f(x_t)
\]

\textbf{Adaptive learning rate:}
\[
E_{t+1} = \beta_2 E_t + (1 - \beta_2) \left( \nabla f(x_t) \right)^2
\]

\textbf{Bias correction:}
\[
\hat{m}_t = \frac{m_t}{1 - \beta_1^t}, \qquad 
\hat{E}_t = \frac{E_t}{1 - \beta_2^t}
\]

\textbf{Parameter update:}
\[
x_{t+1} = x_t - \frac{\hat{m}_t}{\sqrt{\hat{E}_t} + \epsilon}
\]

where \(\beta_1,\beta_2\) are elements of \((0,1)\) typically around \(0.9\).

\subsection*{1.2. Calculus on matrices.}

In machine learning, features of an element of a dataset \(S\) are often represented by an element of a vector space \(\mathbb{R}^p\). If \(S\) contains \(N\) elements, it is represented by an \(N\times p\) matrix \(A\) such that each row of \(A\) represents the features of an element of \(S\). Here are some concepts about calculus in vector spaces and matrices that are important for understanding optimization  in deep learning:

Let $g:\mathbb{R}^n\rightarrow \mathbb{R}^n$ and \(f:\mathbb{R}^n\rightarrow \mathbb{R}\)  differentiable functions defined on $\mathbb{R}^n$. For every element \(x\) of \(\mathbb{R}^n\) and \(u\) in its tangent space, we have:

\begin{align}
\langle\nabla (f\circ g)(x),u\rangle &=d(f\circ g)_x(u)\\
&=df_{g(x)}(dg_x(u))\\
&=\langle\nabla f(g(x)),dg_x(u)\rangle\\
&=\langle dg(x)^T(\nabla_f(g(x)),u\rangle
\end{align}

We deduce that $\nabla (f\circ g)(x)=dg(x)^T(\nabla_f(g(x))$, where $dg(x)^T$ is the transpose of the matrix $dg(x)$.

\bigskip
\bigskip

Let  
\(
M_{m,n}(\mathbb{R}) 
\)
be the space of \( m \times n \) real matrices, and let  
\(
e_{ij} \in M_{m,n}(\mathbb{R})
\)
denote the matrix such that
\[
(e_{ij})_{uv} =
\begin{cases}
1 & \text{if } (u, v) = (i, j), \\
0 & \text{otherwise}.
\end{cases}
\]

 The canonical scalar product in the basis $(e_{ij})_{1\leq i\leq m,1\leq j\leq n}$ coincide with   the scalar product defined by:
 
  \[\langle A,B\rangle=tr(B^TA)
  \]
   where $B^T$ is the transpose of $B$ and $tr$ is the trace. 
 
 This enables to identify the dual space \(M_{m\times n}(\mathbb{R})^*\) of $M_{m\times n}(\mathbb{R})$ with $M_{n\times m}(\mathbb{R})$ by using the isomorphism:

\[
i : M_{n,m}(\mathbb{R}) \rightarrow \left( M_{m,n}(\mathbb{R}) \right)^*
\]
 defined by
\[
i(B)(A) = \operatorname{tr}(BA),
\]

Let \(f:M_{m,n}(\mathbb{R})\rightarrow \mathbb{R}\) be a differentiable function. for every elements $A,B$ of \(M_{m,n}(\mathbb{R})\), $df(A)$ the differential of $f$ at $A$ is an element of \(M_{n\times m}(\mathbb{R})\) and 
 \[df(A)(B)=tr(df(A)B)
 \].
Let \(\nabla f(A)\) be the gradient of $f$ at \(A\), We  also have: 
\[
\langle \nabla f(A),B\rangle=tr(\nabla f(A))^TB)
\]
This implies that $\nabla f(A)$ is $(df(A))^T$.

\section*{2. The Perceptron.}

The perceptron is an algorithm for binary classification. It was invented in 1943 by W. McCulloch and W. Pitts (see [11]). It is a linear classifier: its goal is to determine a hyperplane which separates a finite subset of \( \mathbb{R}^D \) endowed with two classes.

Let \( S \subset B(0,R) \subset \mathbb{R}^D \) be a finite set, and let
\[
f : S \rightarrow \{-1, 1\}
\]
be a labeling function. Suppose that \(\mathbb{R}^D\) is endowed with the canonical scalar product \(\langle,\rangle\).

We say that the pair \((S, f)\) is \textbf{affinely separable} if there exist \( w \in \mathbb{R}^D \) and \( b \in \mathbb{R} \) such that:
\[
\begin{cases}
\langle w, x \rangle + b < 0 & \text{for all } x \in f^{-1}(-1) \\
\langle w, x \rangle + b > 0 & \text{for all } x \in f^{-1}(1)
\end{cases}
\]

That is, the affine hyperplane
\[
H = \left\{ x \in \mathbb{R}^D \,\middle|\, \langle w, x \rangle + b = 0 \right\}
\]
separates the data.

The dataset is \textbf{linearly separable} if \(b=0\).

\subsection*{ Perceptron Algorithm}

\begin{enumerate}
  \item Initialize:
  \[
  w \leftarrow 0 \in \mathbb{R}^D, \quad b \leftarrow 0 \in \mathbb{R}
  \]
  
  \item Repeat:
  \begin{enumerate}
    \item Set \( m \leftarrow 0 \), \(m\) count the misclassified elements of the dataset
    \item  We check if they are misclassified elements. For each \( x \in S \):
    \begin{itemize}
      \item If \( (\langle w, x \rangle + b) f(x) < 0 \), (\(x\) is misclassified) then update:
      \[
      \begin{aligned}
      w &\leftarrow w + f(x)x \\
      b &\leftarrow b + f(x) \\
      m &\leftarrow m + 1
      \end{aligned}
      \]
    \end{itemize}
    \item If \( m = 0 \) after a complete pass through the dataset, all elements are correctly classified and the algorithm terminates. Otherwise, repeat step (b).
  \end{enumerate}
\end{enumerate}

\bigskip

\noindent
 Firstly, we are going to prove that the algorithm terminates after a finite number of steps if the data is linearly separable. Then we are going to extend this result to affinely separable datasets.

\bigskip

 {\bf Theorem 2.1.} {\it Suppose that $(S,f)$  a linearly separable dataset contained in  \(B(0,R)\), the closed ball of radius \(R\), then after a finite number of steps, the perceptron algorithm stops. }
 \medskip

\begin{proof}
We suppose that the separating hyperplane $H$ is defined by:
$\langle u,x\rangle=0$, where $u\in \mathbb{R}^D$, and $\| u\|=1$ is a vector orthogonal to $H$. Let $d$ be the distance between $H$ and $S$. 
We denote by $w_i$ be the weight obtained after $i$ updates.
Suppose that the element $x_i\in S$ is a $i$-misclassified element of $S$, we update the weight to $w_i=w_{i-1}+f(x_i)x_i$.

We have:

\begin{align}
\langle w_i,u\rangle &=\langle w_{i-1} +f(x_i)x_i,u\rangle\\
                     &=\langle w_{i-1},u\rangle+f(x_i)\langle x_i,u\rangle 
\end{align}

We can write $x_i=x_i^0+x_i^1$ where $x_i^0=tu$ and $x_i^1$ is an element of $H$. Since $d$ is the distance between $H$ and $S$, we have:

 $f(x_i)\langle x_i,u\rangle=f(x_i)\langle tu,u\rangle=\mid t\mid=d(x_i,H)\geq d$,
 since we have assumed, $\| u\|=1$.  We deduce that:

\begin{align}
\langle w_i,u\rangle &=\langle w_{i-1} +f(x_i)x_i,u\rangle\\
                     &\geq \langle w_{i-1},u\rangle+d. 
\end{align}

We also have:

\begin{align}
\langle w_i,w_i\rangle=&\langle w_{i-1}+f(x_i)x_i,w_{i-1}+f(x_i)x_i\rangle\\
                       &=\langle w_{i-1},w_{i-1}\rangle +2f(x_i)\langle w_{i-1},x_i\rangle +\langle x_i,x_i\rangle\\
                       &\leq \langle w_{i-1},w_{i-1}\rangle+R^2 
\end{align}

This follows from the facts that $f(x_i)\langle w_{i-1},x_i\rangle\leq 0$ since $x_i$ is misclassified, and $\langle x_i,x_i\rangle\leq R^2$ since $x_i\in B(0,R)$.

This applied after $M$ updates,

$\langle w_M,u\rangle\geq Md$

$\| w_M\|^2\leq R^2M$.

Cauchy-Schwarcz implies that $\mid \langle w_M,u\rangle\ \mid \leq \|w_M\|\|u\|$. We deduce that:
  $Md\leq \| w_M\|\leq \sqrt{R^2M}$.

This implies that $M\leq { R^2\over d^2}$. We deduce that we cannot repeat the algorithm more than ${ R^2\over d^2}$ times.

\end{proof}

{\bf Corollary 2.2.}
If $(S,f)$ is affinely separable, after a finite number of steps, the perceptron algorithm stops. 
\medskip

\begin{proof}
We embed the dataset \(S\) to \(\mathbb{R}^{D+1}\) using the map:

$i:\mathbb{R}^{D}\longrightarrow \mathbb{R}^{D+1}$

$i(x)=(x,1)$.

Suppose that \(S\) is affinely separable, and the separating hyperplane is defined by \(w\in\mathbb{R}^D\) and \(b\in \mathbb{R}\), then \(i(S)\) is linearly separable, and its separating hyperplane is defined \((w,b)\).
The result follows from the fact that
the \(i\)-weight \(\hat w_i\) obtained by applying the perceptron algorithm to \(i(S)\) is 
\((w_i,b_i)\), where \(w_i\) and \(b_i\) are the \(i\)-weights obtained by applying the perceptron algorithm to \(S\). 

\end{proof}

{\bf Remark.}

\medskip
While proving the previous theorem, we have seen that the maximal number of updating steps is proportional to \( R^2\), where \(R\) is an upper boun of the norm of elements of the dataset. This shows that normalizing the dataset may increase performance. 
 In practice the infinity norm is used to normalize the dataset: each feature is normalized individually to make them directly comparable. This is achieved with:

\medskip

{\bf Min Max scaler.}

Suppose that the dataset is represented by a \(N\times D\) matrix \(A\).  We replace the column $x$ of \(A\) by $x'$ defined by the formula:

$$
x'={{x-x_{min}}\over {x_{max}-x_{min}}}
$$
where $x_{min}$ and $x_{max}$ are respectively the minimum and the maximum value of the feature column.

Another algorithm used in machine learning to normalize datasets is:

\medskip

{\bf Standard scaler.}
$$
x'={{x-\mu}\over\sigma}
$$

where $\mu$ is the mean of the feature $x$ and $\sigma$ its standard deviation of the column \(x\). It rescales the feature $x$ to a feature $x'$ whose mean is $0$ and whose standard deviation is $1$. 
\medskip

\section*{3. Logistic Regression.}

A high percentage, but not all, of a labelled dataset \((S,f)\) may be separable by a hyperplane \( H \). In such cases, the perceptron algorithm will not converge. One modification to improve performance is to introduce a loop that repeats the algorithm a fixed number of times. 
However, a more effective approach is \emph{logistic regression}, which handles non-separable data more efficiently by optimizing a  loss function and providing probabilistic outputs for classification.

We suppose that the dataset $S$ is a subset  of ${\bf R}^D$ which contains $N$ elements. Equivalently, it can be represented with a $N\times D$ matrix whose rows are elements of $S$. 
For convenience, we are going to suppose here that $f$ takes its values in \(\{0,1\}\).
\bigskip

For a given hyperplane \( H \subset \mathbb{R}^D \), defined by the equation \( l_H(x) = 0 \), where \( l_H(x) \) is an affine function, we define the forward map associated to \(H\) to be the restriction  \( f_H(x)\) of  \(l_H(x) \)  to the dataset \( S \).

Let the sigmoid function \( \sigma : \mathbb{R} \to (0, 1) \) be defined as:
\[
\sigma(z) = \frac{1}{1 + e^{-z}}
\]

Then, the contribution of a point \( x \in S \) to the loss function is:
\[
L_H(x) = -f(x) \log\left(\sigma(f_H(x))\right) - \left(1 - f(x)\right) \log\left(1 - \sigma(f_H(x))\right)
\]

The total logistic loss over the dataset is:
\[
L_H = \frac{1}{N} \sum_{x \in S} L_H(x)
\]

This loss function allows us to evaluate how well the hyperplane \( H \) classifies the elements of \( S \).

The term \( \sigma(f_H(x)) \) represents the model's estimated probability that the label \( f(x) \) is 1.  
 Suppose \( f(x) = 1 \). Then the contribution of \(x\) to the loss becomes:
\[
L_H(x) = -\log(\sigma(f_H(x)))
\]
If \( \sigma(f_H(x)) \) is close to 1, the loss is close to 0, indicating a correct and confident prediction. If \( \sigma(f_H(x)) \) is close to 0, the loss becomes large, indicating a confident but incorrect prediction.

On the other hand, \( 1 - \sigma(f_H(x)) \) represents the estimated probability that \( f(x) = 0 \).  
 Suppose \( f(x) = 0 \). Then the contribution of \(x\) to the loss becomes:
\[
L_H(x) = -\log(1 - \sigma(f_H(x)))
\]
If \( \sigma(f_H(x)) \) is close to \(0\), the loss is again small. But if \( \sigma(f_H(x)) \) is close to 1, the loss becomes large.

Thus, the loss is small when the model makes confident correct predictions, and large when it is confidently wrong.

We can write 
\[
f_H(x_1,\ldots,x_D) = w_1x_1 + \ldots + w_Dx_D + b,
\]
and observe that \( L_H \) is a differentiable function of \( w_1,\ldots,w_D, b \). The hyperplane \( L_H \) gives a better classification than \( L_{H'} \) if \( L_H < L_{H'} \). The global minimum, if it exists, of the function \( L_H(w_1,\ldots,w_D,b) \), defines a hyperplane \( H_0 \) which is the best classifier of \( (S,f) \).

In practice, it is not always possible to compute the exact equation of \( H_0 \), but with gradient descent (See section 1.1), we can find hyperplanes that classify a high percentage of \( (S,f) \) if it exists. For this purpose, we compute the differential of \( L_H \).

Let \( W = (w_1,\ldots,w_D) \in (\mathbb{R}^D)^* \), and \( b \in \mathbb{R} \), to which we associate \( b_N \in \mathbb{R}^N \) with all coordinates equal to \( b \). Let \( M_{N\times D}(\mathbb{R}) \) be the space of \( N \times D \) matrices. Define the map
\[
F_{W,b}: M_{N\times D}(\mathbb{R}) \to \mathbb{R}^N, \quad F_H(X) = XW + b_N.
\]
Define also \( \sigma_N: \mathbb{R}^N \to \mathbb{R}^N \) by
\[
\sigma_N(x_1,\ldots,x_N) = (\sigma(x_1),\ldots,\sigma(x_N)),
\]
and
\[
L_Y: \mathbb{R}^N \to \mathbb{R}, \quad L_Y(h_1,\ldots,h_N) = -\frac{1}{N} \sum_{i=1}^N \left[ y_i \log(h_i) + (1 - y_i)\log(1 - h_i) \right].
\]
 Where \(y_i=f(x_i)\). We have:
\[
L_H = L_Y \circ \sigma_N \circ F_{W,b}.
\]
To apply gradient descent, we compute the gradient of \( L_H \).

Since \( F_{W,b} \) is affine in \( (W,b) \), for every \(U\in \mathbb{R}^D\) and \(c\in \mathbb{R}\), we have:
\[
\frac{\partial F_{W,b}}{\partial W}(U) = XU, \quad \frac{\partial F_{W,b}}{\partial b}(c) = c \cdot (1,\ldots,1) \in \mathbb{R}^N.
\]

The Jacobian matrix of \( \sigma_N \) is:
\[
d(\sigma_N) = \begin{bmatrix}
\sigma(x_1)(1 - \sigma(x_1)) & 0 & \cdots & 0 \\
0 & \sigma(x_2)(1 - \sigma(x_2)) & \cdots & 0 \\
\vdots & \vdots & \ddots & \vdots \\
0 & 0 & \cdots & \sigma(x_n)(1 - \sigma(x_n))
\end{bmatrix}.
\]

The differential of \( L_Y \) is the one-form:
\[
dL_Y = -\frac{1}{N} \sum_{i=1}^N \left( \frac{y_i}{h_i} - \frac{1 - y_i}{1 - h_i} \right) dh_i 
= -\frac{1}{N} \sum_{i=1}^N \frac{y_i - h_i}{h_i(1 - h_i)} dh_i.
\]

Setting \( h_i = \sigma(x_i) \), we deduce:
\[
dL_Y \circ d\sigma_N = \frac{1}{N} \sum_{i=1}^N (\sigma(x_i) - y_i) dx_i.
\]

Let \( \hat{Y} = (\sigma(f_H(x_1)), \ldots, \sigma(f_H(x_n))) \), The section 1.2 implies that: 
\[
dL_Y \circ \sigma_N = \frac{1}{N} (\hat{Y} - Y)^T.
\]

Thus, the differential of \( L_H \) in direction \( U \) is:
\[
dL_H(U) = \frac{1}{N} (\hat{Y} - Y)^T XU = {1\over N}\operatorname{tr}\left((X^T(\hat{Y} - Y))^T U\right).
\]

Similarly, we have:
\[
dL_H(b) = {1\over N}\sum_{i=1}^N (\hat{Y}_i - Y_i).
\]

The section 1.2 shows that the gradient \(L_H\) of \(L_H\) is: 
\[
\nabla L_H(W, b) = \left( \nabla_W L_H, \nabla_b L_H \right) 
= \left({1\over N} X^T(\hat{Y} - Y), {1\over N}\sum_{i=1}^N (\hat{Y}_i - Y_i) \right).
\]

\bigskip

\textbf{Logistic regression algorithm.}

\medskip

We present the logistic regression algorithm.

\medskip

\textbf{Initializing weights:}

We initialize the weight vector \( W \) using a normal distribution, and set \( b \) to zero.

\medskip

\textbf{Feed forward:}

We define a function to compute \( \hat{Y} \), the output of the model. It is similar to the function used in the perceptron algorithm:

\[
F_{W,b}(X) = \sigma_N(XW + b_N)
\]

\medskip

\textbf{Training loop:}

We choose the number of epochs \( n_{\text{epochs}} \), and a learning rate \( \alpha \) to update the weights.  
In Python, we can create an empty list \texttt{losses} to keep track of the successive losses.

We repeat the following steps \( n_{\text{epochs}} \) times:

\begin{itemize}
    \item Compute the prediction: \( \hat{Y} = F_{W,b}(X) \)
    \item Compute the loss: \( \text{loss} = L_H(\hat{Y}, Y) \)
    \item Append the loss to \texttt{losses}
    \item Update the weights:
    \[
    W \leftarrow W - \alpha \cdot \frac{1}{N} X^T (\hat{Y} - Y)
    \]
    \[
    b \leftarrow b - \alpha \cdot \frac{1}{N} \sum_{i=1}^N (\hat{Y}_i - Y_i)
    \]
\end{itemize}

\medskip

After the training loop, we can plot the list of losses to visualize whether it decreases, that is, whether the model learns.

\bigskip

There exist datasets \( (S, f) \) that are separable but not linearly separable. For example, consider \( S \) as the union of:

\begin{itemize}
    \item 100 points in \( \mathbb{R}^2 \), uniformly distributed in the ball of radius 1, labeled \( 0 \),
    \item 100 points uniformly distributed in the annulus between circles of radius 1 and 2, labeled \( 1 \).
\end{itemize}

A possible solution to classify this dataset is to embed it into a higher-dimensional vector space where they become linearly separable. This can be  achieved because  we have prior information about the topology of the dataset that allows us to construct such an embedding.

However, with real-world data such as images, where features may lie in a space of dimension \( 1024 \times 1024 \) this approach is not feasible.
Neural networks with more than one hidden layer are able to classify datasets that are not linearly separable.

\section*{ 4. Deep Neural Networks.}

\textbf{Definition 4.1.}  
A \emph{neural network of depth} \( d+1 \) is defined by a finite set \( G \), together with a partition  into layers:
\[
G_0, G_1, \dots, G_d.
\]

Let \( A \) be such a neural network of depth \( d+1 \). We denote by \( n_i = |G_i| \) the number of nodes in the \( i \)-th layer.

\medskip

Let \( A = (G_0, G_1, \dots, G_d) \) be a neural network. Consider a dataset \( S \subset \mathbb{R}^{n_0} \) consisting of \( N \) elements. Suppose there exists a \emph{label map}:
\[
Y : S \rightarrow \{0, 1, \dots, n_d-1\},
\]
 The goal is to construct, using the neural network \( A \), a \emph{prediction map}:
\[
\hat{Y} : S \rightarrow \mathbb{R}^{n_d},
\]
such that for every \( x \in S \), the \( i \)-th coordinate of \( \hat{Y}(x) \) represents the predicted probability that the label of \( x \) is \( i \).

Typically, the labeled dataset \( X \) is a finite subset of a larger dataset \( \hat{X} \), which is defined conceptually and whose full content is not known. For example, consider a dataset consisting of images of cats and dogs. Each image can be represented as three matrices corresponding to the red, green, and blue (RGB) channels, allowing us to embed each image in a high-dimensional real vector space.
Clearly, at any fixed time \( t \), the dataset \( X \) cannot contain all images of cats that have ever existed, exist at time \( t \), or will exist in the future. We define the label map \( Y \) so that it assigns \( 0 \) to cats and \( 1 \) to dogs. We expect the prediction map \( \hat{Y} \) to generalize to unseen examples. That is, for any new cat image represented as a point \( c \in \mathbb{R}^{n_0} \), we expect
\(
\hat{Y}(c) \) to be close to \(0\).

To define $\hat Y$, we firstly define the following functions:

\paragraph{The ReLU.}
We denote by $r : \mathbb{R} \to \mathbb{R}$ the function defined by
\[
r(x) = \begin{cases}
x & \text{if } x \geq 0, \\
0 & \text{if } x < 0.
\end{cases}
\]
This function is commonly known as the \emph{rectified linear unit}, or \emph{ReLU}.

\medskip

\paragraph{The softmax.}

\medskip

Let \( s : \mathbb{R}^{n_d} \to \mathbb{R}^{n_d} \) be the function defined by
\[
s(x_0, \dots, x_{d-1}) = \left( \frac{e^{x_0}}{\sum_{i=0}^{d-1} e^{x_i}}, \dots, \frac{e^{x_{d-1}}}{\sum_{i=0}^{d-1} e^{x_i}} \right).
\]

\subsection*{Recursive Definition of a Forward map.}

Let \( W_i \) be a real \(n_{i-1}\times n_i\) (weight) matrix and \( b_i \in \mathbb{R}^{n_i} \) a bias vector for the layer \( i \), \( i\in \{1,\cdots,d \}\). Define the forward map recursively as follows:

\paragraph{Initialization}
\[
H_0 = X,
\]
where \( X \) is the \( N \times n_0 \)  input matrix, with \( N \) samples and \( n_0 \) features.

\paragraph{First Layer}
\[
Z_1 = H_0 W_1 + b_1,
\]
\[
H_1 = a_1(Z_1),
\]
where the activation function \( a_1 \) acts componentwise:
\[
a_1 = 
\begin{cases}
\text{ReLU}, & \text{if } 1 < d, \\
\text{softmax}, & \text{if } 1 = d.
\end{cases}
\]

\paragraph{Recursive Step}
Suppose \( H_i \) and \( Z_i \) are defined for \( i < d \). Then for \( i+1 \leq d \), define:
\[
Z_{i+1} = H_i W_{i+1} + b_{i+1},
\]
\[
H_{i+1} = a_{i+1}(Z_{i+1}),
\]
where the activation function \( a_{i+1} \) acts componentwise and is defined as:
\[
a_{i+1} = 
\begin{cases}
\text{ReLU}, & \text{if } i+1 < d, \\
\text{softmax}, & \text{if } i+1 = d.
\end{cases}
\]

\paragraph{Output and Prediction}

Let \( \hat{Y} = H_d \in \mathbb{R}^{N \times n_d} \) denote the final output of the network, where \( d \) is the total number of layers. Each row of \( \hat{Y} \) corresponds to a probability distribution over \( n_d \) classes, i.e., the components of each row sum to 1 due to the softmax activation function at the final layer:
\[
\sum_{j=1}^{n_d} \hat{Y}_{l j} = 1 \quad \text{for all } l = 1, \dots, N.
\]

If \( x_l \) is the \( l \)-th input row of \( X \), then the predicted class label \( \hat{y}_l \) for \( x_l \) is defined by:
\[
\hat{y}_l = \arg\max_{j = 1, \dots, n_d} \hat{Y}_{l j}.
\]

This selects the class with the highest predicted probability for each input sample.

\medskip

We have defined a function:

\[
F_{W, b} : M(n_0, n_1)(\mathbb{R}) \times \cdots \times M(n_{d-1}, n_d)(\mathbb{R}) \times \mathbb{R}^{n_1} \times \cdots \times \mathbb{R}^{n_d} \longrightarrow \mathbb{R}^{n_d}
\]

such that

\[
F_{W, b}(W_1, \ldots, W_{n_d}, b_1, \ldots, b_{n_d}) = \hat{Y}
\]

We are going to define a loss function which enables to quantify the accuracy of the predictions given \(W_1,\cdots,W_d,b_1,\cdots b_{n_d}\).

We assume that labels are one hot encoded equivalently, $i$ corresponds to the element \(e_i\) of \( \mathbb{R}^{n_d}$ whose only non zero coordinate is its $i$-coordinate which is equal to $1$.
The label of the element $x_i$ of $X$ will be written \( y_i=(y_{i1}\cdots y_{in_d})\).

\paragraph{Cross Entropy Loss}

The prediction \( \hat{Y} \) is an \( N \times n_d \) matrix, where the \( i \)-th row (denoted \( \hat{y}_i = (\hat{y}_{i1}, \cdots, \hat{y}_{in_d}) \)) represents the prediction for the input element \( x_i \in X \).

The loss function we use is the cross-entropy loss, defined as:

\[
L(\hat{Y}, Y) = -\frac{1}{N} \sum_{i=1}^{N} \sum_{j=1}^{n_d} y_{ij} \log(\hat{y}_{ij})
\]

Suppose that the label of the element \(x_i\) of the dataset is \(e_p\) its contribution to the 
loss is \( -\frac{1}{N} log(\hat y_{ip})$. If \( \hat y_{ip}\) is close to \(1\), its contribution to the loss is very small. It is natural to find the values of the parameters \( W_1,\cdots,W_{n_d},b_1,\cdots b_{n_d}$ for which the loss is minimum. We are going to approach these values with backpropagation and gradient descent.

\subsection*{Back Propagation and Gradient Descent.}

Backpropagation is the process of computing the derivative and the gradient of the loss function relatively to
\( W_1,\cdots W_{d},b_1,\cdots,b_{d}\) by starting with \(W_{n_d}\) and \(b_{d}\) and going backward.

The composition of the Jacobian matrices of derivatives gives:

\[
\frac{d L}{d W_{d}} = \frac{d L}{d H_{d}} \circ \frac{d H_{d}}{d Z_{d}} \circ \frac{d Z_{d}}{d W_{d}}
\]

\( \frac{d L}{d H_{d}}\) is the $1$-form:
\[
 -\frac{1}{N}\sum_{i=1}^{i=N}\sum_{j=1}^{j=n_d}\frac{y_{ij}}{\hat y_{ij}}d\hat y_{ij}
 \]
 To compute \( \frac{d H_{d}}{d Z_{d}}\), We need the following lemma:

\bigskip 
{\bf Lemma 4.1.}
\medskip

Let \( s:\mathbb{R}^m\longrightarrow \mathbb{R}^m\) the softmax,
Write \( s(x_1,\cdots ,x_m)=(s_1,\cdots ,s_m)\). The differential of \(s\) is:

is:\[
 \begin{pmatrix}
s_1(1-s_1) & -s_1s_2 & \cdots & -s_1s_m \\
-s_1s_2 & s_2(1-s_2) & \cdots & -s_2s_m\\
\vdots & \ddots & & \vdots \\
\vdots & & \ddots & \vdots \\
-s_1s_m & \cdots &  & s_m(1-s_m)
\end{pmatrix}
\]

\medskip

We have:

\begin{proof}

We have:
\begin{align}
\frac{\partial s_i}{\partial x_i} &=\frac{e^{x_i}(e^{x_1}+\cdots +e^{x_m})-(e^{x_i})^2}{(e^{x_1}\cdots e^{x_m})^2}\\
&=\frac{e^{x_i}}{e^{x_1}+\cdots+ e^{x_m}}\frac{e^{x_1}+\cdots +e^{x_m}-e^{x_i}}{e^{x_1}+\cdots + e^{x_m}}\\
&=s_i(1-s_i)
\end{align}

If \(i\neq j\)

\begin{align}
\frac{\partial s_i}{\partial x_j} &= \frac{-e^{x_i}e^{x_j}}{(e^{x_1}+\cdots+e^{x_m})^2}\\
&=-s_is_j
\end{align}
\end{proof}

We deduce that \( \frac{dH_{d}}{dZ_{d}}\) is the \( N\times n_d\times N\times n_d\) matrix:

\[
\begin{pmatrix}
M_1 & 0 & \cdots & 0 \\
\vdots & \ddots & & \vdots \\
\vdots & & \ddots & \vdots \\
0 & \cdots & 0 & M_N
\end{pmatrix}
\]

where

\[
 M_i=\begin{pmatrix}
\hat y_{i1}(1-\hat y_{i1}) & -\hat y_{i1}\hat y_{i2} & \cdots & -\hat y_{i1}\hat y_{in_d} \\
-\hat y_{i1}\hat y_{i2} & \hat y_{i2}(1-\hat y_{i2}) & \cdots & -\hat y_{i2}\hat y_{in_d}\\
\vdots & \ddots & & \vdots \\
\vdots & & \ddots & \vdots \\
-\hat y_{i1}\hat y_{in_d} & \cdots &  & \hat y_{in_d}(1-\hat y_{in_d})
\end{pmatrix}
\]

We deduce that \( \frac{dL}{dH_{n_d}}\circ \frac{dH_{n_d}}{dZ_{n_d}}\) is the \(1\)-form defines on the space of \(N\times n_d\) matrices, that is an \(n_d\times N\) matrix whose coefficient \(i,j\) is

\begin{align}
&-\frac{1}{N}[-\hat y_{i1}\hat y_{ij}\frac{y_{i1}}{\hat y_{i1}}-\cdots+\hat y_{ij}(1-\hat y_{ij})\frac{y_{ij}}{\hat y_{ij}}-\cdots -\hat y_{ij}\hat y_{in_d}\frac{y_{in_d}}{\hat y_{n_d}}]\\
&=\frac{1}{N}(\hat y_{ij}-y_{ij})
\end{align}
\bigskip

since \(\sum_{j=1}^{j=n_d}y_{ij}=1\)

We deduce that \( \frac{dL}{dH_{n_d}}\circ \frac{dH_{n_d}}{dZ_{n_d}}=\frac{1}{N}(\hat Y-Y)^T\).

Since \( Z_{n_d}=H_{n_{d}-1}W_{n_d}+b_{n_d}$, we deduce that \( \frac{dZ_{n_d}}{dW_{n_d}}\) is the linear map defined on the space on \( N\times n_d\) matrices by:

\[
\frac{dZ_{n_d}}{dW_{n_d}}(U)=H_{n_d-1}U
\]

This implies that:

\begin{align}
  \frac{dL}{dW_{n_d}}(U) &=\frac{dL}{dH_{n_d}}\circ \frac{dH_{n_d}}{dZ_{n_d}}\circ \frac{dZ_{n_d}}{dW_{n_d}}(U)\\
 &=tr((\hat Y-Y)^TH_{n_d-1}U)\\
 &=tr((H_{n_d-1}^T(\hat Y-Y))^TU)
\end{align}
 
 We deduce that \( \nabla L({W_{n_d}})=H_{n_d-1}^T(\hat Y-Y)\)

We also have \( \frac{dL}{db_{n_d}}=\frac{dL}{dH_{n_d}}\circ \frac{dH_{n_d}}{dZ_{n_d}}\circ \frac{dZ_{n_d}}{db_{n_d}}\)

For every element \(v=(v_1,\cdots v_{n_d})\) of \(\mathbb{R}^{n_d}\), \( \frac{dZ_{n_d}}{db_{n_d}}(v)\) is the
\(N\times n_d\) \(v_N\) matrix whose \(i\)-row is \(v\). We deduce that:

\begin{align}
\frac{dL}{db_{n_d}}&=tr((\hat Y-Y)^Tv_N)\\
&=(\sum_{i=1}^{i=N}(\hat Y-Y)_{i1}v_1,\cdots ,\sum_{i=1}^{i=N}(\hat Y-Y)_{in_d}v_{n_d})
\end{align}

This is equivalent to saying that \( \nabla L_{b_d}=\sum_{i=1}^{i=N}(\hat Y-Y,axis=1)\) which is the \(\mathbb{R}^{n_d}\) vector obtained by summing  the \(N\times n_d\) matrix \(\hat Y-Y\) along its columns.

\bigskip

Suppose known \(\frac{dL}{dH_{i+1}}, \frac{dL}{dZ_{i+1}},\frac{dL}{dW_{i+1}}, \frac{dL}{db_{i+1}}\), since \( Z_{i+1}=H_iW_{i+1}+b_{i+1}\), we deduce that for every vector of \( M_{N\times n_i}(\mathbb{R})\), we have:

\begin{align}
\frac{dL}{dH_i}(U)&=\frac{dL}{dZ_{i+1}}( \frac{dZ_{i+1}}{dH_i}(U))\\
                  &=\frac{dL}{dZ_{i+1}}(UW_{i+1})\\
                  &=tr(\nabla L ({Z_{i+1}}))^TUW_{i+1})\\
                  &=tr((\nabla L ({Z_{i+1}})W_{i+1}^T)^TU)
\end{align}

We deduce that \(\nabla L({H_i})=\nabla L({Z_{i+1}})W_{i+1}^T\)

To compute \(\nabla L(Z_i)\) we need to compute the derivative of the ReLU. We are going to consider it has a distribution defined on test functions with finite support.

\bigskip

{\bf Proposition 5.2.}
The derivative of the ReLu is the function defined by:

\[
r'(x) = \begin{cases}
1 & \text{if } x \geq 0, \\
0 & \text{if } x < 0.
\end{cases}
\]

\begin{proof}
Let \(\phi\) a test function, the derivative of  the distribution \(r\) is defined by:
\begin{align}
\int_{\infty}^{+\infty}r'\phi
                              &=-\int_{\infty}^{+\infty}r\phi'\\ 
                              &=-\int_0^{+\infty}x\phi'\\
                               &=[-x\phi]_0^{+\infty}+\int_0^{+\infty}\phi\\
                               &=\int_0^{+\infty}\phi
\end{align}

We deduce that the derivative of the ReLU is the distribution \(r'\).
\end{proof}

We have \( \frac{dL}{dZ_i}=\frac{dL}{dH_i}\circ \frac{dH_i}{dZ_i}\)

The derivative \(\frac{dH_i}{dZ_i}\) is the \(N\times n_i\times N\times n_i\) diagonal matrix whose entry at position \((jN+k,jN+k)\) is \(r'(z_{jk})\), where \(z_{jk}\) is the entries of \(Z_i\) at position \((j,k)\). This implies that for every element \(U\) of \(M_{N\times n_i}(\mathbb{R})\), 
\[\frac{dH_i}{dZ_i}(U)=r'(Z_i)*U
\] where \(r'(Z_i)\) is the \(N\times n_i\) matrix whose entries are \(r'(z_{jk})\) and \(r'(Z_i)*U\) is the Hadamard product: that is the matrix obtained by multiplying \(r'(Z_i)\) and \(U\) coefficientwise.

 We deduce that:

\begin{align}
\frac{dL}{dZ_i}(U)&=\frac{dL}{dH_i}(r'(Z_i)*U)\\
                  &=tr(\nabla L(Z_i)^Tr'(Z_i)*U)\\
                  &=tr((\nabla L(Z_i)*r'(Z_i))^TU)\\
\end{align}

We deduce that:

\[
\nabla L(Z_i)=\nabla L(H_i)*r'(Z_i)
\]

We have \(Z_i=H_{i-1}W_i+b_i\). This implies that \(\frac{dZ_i}{dW_i} \) is the linear map of the space of \(N\times n_i\) matrices defined by: \(\frac{dZ_i}{dW_i}(U)=H_{i-1}U\). We deduce that:

\begin{align}
\frac{dL}{dW_i}(U) &= \frac{dL}{dZ_i}(H_{i-1}U)\\
                   &=tr((\nabla L(Z_i)^TH_{i-1}U)\\
                   &=tr((H_{i-1}^T\nabla L(Z_i))^TU)
\end{align}

We deduce that \(\nabla L(W_i)=H_{i-1}^T\nabla L(Z_i)\).

\medskip

We also have \( \frac{dL}{db_{i}}=\frac{dL}{dZ_i}\circ \frac{dZ_{i}}{db_{i}}\)

For every element \(v=(v_1,\cdots v_{n_i})\) of \(\mathbb{R}^{n_i}\), \( \frac{dZ_{i}}{db_{i}}(v)\) is the
\(N\times n_i\) matrix \(v_N\) whose \(j\)-row is \(v\). We deduce that:

\begin{align}
\frac{dL}{db_{i}}(v)&=tr((\nabla L(Z_i))^Tv_N)\\
&=(\sum_{j=1}^{j=N}(\nabla L(Z_i))_{j1}v_1,\cdots, \sum_{j=1}^{j=N}(\nabla L(Z_i))_{jn_i}v_{n_i})
\end{align}

This is equivalent to saying that \( \nabla L({b_i})=\sum_{j=1}^{j=N}(\nabla L(Z_i),axis=1)\)

\subsection*{Deep Neural Network Algorithm.}

{\bf Weight Initialization}

\medskip

In order to implement a deep neural network, we must firstly initialize the weights, which consist of matrices \(W_i\) and bias vectors \(b_i\).
There are several methods to initialize the matrices \(W_i\). One common approach is to sample entries from a standard normal distribution using \texttt{np.random.randn} in NumPy.
The bias vectors \(b_i\) are typically initialized to zero.

\medskip

\noindent\textbf{Remark.}

\medskip

We do not initialize the weights \(W_i\) as zero matrices, since this would imply that every neuron receives the same input. If \(W_1\) and \(b_1\) are zero vectors
 all coordinates of
\(H_1 = a_1(Z_1)\)
are equal. Recursively, we deduce that the coordinates of \(H_i\) remain constant if all \(W_1, \ldots, W_i\) and \(b_1, \ldots, b_i\)
 are initialized to zeros. As a result, every neuron in a given layer produces the same output,  and no meaningful learning can occur if the network has more than two hidden layers if the depth of the layer is strictly superior to \(1\) 
 since  \(\nabla L({H_i})=\nabla L({Z_{i+1}})W_{i+1}^T\) and \(W_{i+1}=0$ implies that \(\nabla L(H_i)=0\). This failure to break symmetry prevents the network from converging to a useful solution.

\medskip

{\bf Forward Map}

\medskip
Using the recursive relations:
\[
Z_i = H_{i-1} W_i + b_i, \quad H_i = a_i(Z_i),
\]
we define a function that stores the intermediate values \(Z_i\) and \(H_i\) into two lists, denoted by \(Z\) and \(H\) since the values of \(H\)  values are essential for backpropagation.

\medskip

{\bf The Training Loop}

\medskip

We repeat the following steps for a fixed number of iterations:
\begin{itemize}
  \item Compute the forward map;
  \item Compute the loss.

 \item Backpropagation and Weight Update: compute the gradients and update the weights accordingly.
\end{itemize}
Finally, we compute the accuracy to evaluate the percentage of inputs correctly classified by the network.

\subsection*{5. Universal Approximation Theorem and Overfitting.}

A neural network with \(n\) input nodes and \(m\) output nodes  defines a function from \(\mathbb{R}^n\) to \( \mathbb{R}^m\). The universal approximation theorem shows that every function defined on a compact subset of \(\mathbb{R}^n\) can be approximated arbitrarily well by a neural network with one hidden layer and an activation function \(\sigma\) which satisfies some boundary conditions. Recall that \(\sigma\) is a function whose domain and codomain is $\mathbb{R}$. If there is no confusion, we also denote  by \(\sigma\) the function which acts on \(\mathbb{R}^p\) by \(\sigma(x_1,\cdots,x_p)=(\sigma(x_1),\cdots,\sigma(x_p))\).

\bigskip

{\bf Theorem 5.1.}(Cybenko [4])
Let \( \sigma: \mathbb{R} \to \mathbb{R} \) be a continuous function satisfying
\[
\lim_{t \to -\infty} \sigma(t) = 0, \quad \lim_{t \to +\infty} \sigma(t) = 1.
\]
Let \( f \) be a continuous function defined on a compact subset \( K \subset \mathbb{R}^n \) which takes it values in \(\mathbb{R}^m\). Then for every \( \epsilon > 0 \), there exist:

\begin{itemize}
    \item an integer \( l \in \mathbb{N} \),
    \item a  matrix \( A \in M_{l\times n}(\mathbb{R})\),
    \item vectors \( b \in \mathbb{R}^l \),
    \item a matrix \(C\in M_{m\times l}(\mathbb{R})\), for \( i = 1, \dots, m \),
\end{itemize}

such that the function
\[
g(x) =  C \sigma(A x + b)
\]
satisfies
\[
\sup_{x \in K} \|f(x) - g(x)\| < \epsilon.
\]

\medskip

\begin{proof}

It is enough to show the theorem for \(m=1\).

Let \(a_i\) the \(i\)-row of a matrix \(A\). We aim to show that the set of finite linear combinations of the form
\[
\mathcal{N}_\sigma := \left\{ x \mapsto \sum_{i=1}^{i=l}c_i\sigma(a_i^Tx+b_i) \right\}
\]
is dense in the space of continuous functions defined on \(K\), \( C(K) \),
where \(c_i,b_i\in \mathbb{R}\) and \(a_i\in\mathbb{R}^n\).

Suppose \( \mathcal{N}_\sigma \) is not dense in \( C(K) \). Then, by the Hahn–Banach theorem, there exists a nonzero bounded linear functional \( L \in C(K)^* \) such that
\[
L(g) = 0 \quad \text{for all } g \in \mathcal{N}_\sigma.
\]

By the Riesz Representation Theorem, \( L \) corresponds to a signed Borel measure \( \mu \) such that
\[
L(f) = \int_K f(x) \, d\mu(x).
\]
Thus, we have
\[
\int_K \sigma(a^T x + b) \, d\mu(x) = 0 \quad \text{for all } a \in \mathbb{R}^n, b \in \mathbb{R}.
\]

Let $H = \{x \in \mathbb{R}^n : a^T x + b > 0\}$ be an open half-space. Define for $u > 0$ the rescaled function
\[
\sigma_u(t) := \sigma(ut),
\]
so that
\[
\sigma_u(a^T x + b) \to \mathbf{1}_H(x) \quad \text{pointwise as } u \to +\infty.
\]
Since $\sigma$ is bounded,   the functions $\sigma_u(a^T x + b)$ are uniformly bounded. By the assumption,
\[
\int_K \sigma(u(a^T x + b))\, d\mu(x) = 0 \quad \text{for all } u > 0.
\]
By the bounded convergence theorem,
\[
\mu(H \cap K) = \lim_{u \to +\infty} \int_K \sigma_u(a^T x + b)\, d\mu(x) = 0.
\]
Thus, $\mu$ assigns zero measure to all intersections of $K$ with open half-spaces.

 Define the collection
\[
\mathcal{H} := \{ H \cap K : H \text{ is a half-space in } \mathbb{R}^n \}.
\]
This is a $\pi$-system that generates the Borel $\sigma$-algebra on $K$. Since $\mu$ is a finite signed measure and vanishes on this generating $\pi$-system, we deduce that:
\[ 
\mu(E) = 0 \quad \text{for all Borel sets } E \subseteq K.
\]

Therefore, $\mu = 0$, contradiction.

\end{proof}

\bigskip

{\bf Corollary 5.2.}
 Let \( K \subset \mathbb{R}^n \) be a compact set. The vector subspace \(V(RELU)\) of \( C(K, \mathbb{R}^m) \), the space of continuous functions from \(K\) to \(\mathbb{R}^m\), generated by the family
\[
C \, \mathrm{ReLU}(A x + b),
\]
where \( A \in M_{l\times n}(\mathbb{R}) \), \( C \in M_{m\times l}(\mathbb{R}) \), and \( b \in \mathbb{R}^l \) is dense in \( C(K, \mathbb{R}^m) \).

\medskip

\begin{proof}
Let \(\sigma(z) = \mathrm{ReLU}(z) - \mathrm{ReLU}(z - 1)\). Explicitly, this function is defined as
\[
\sigma(z) = 
\begin{cases}
0, & z < 0, \\
z, & 0 \leq z < 1, \\
1, & z \geq 1.
\end{cases}
\]

The theorem 5.1 implies that the vector space \(V(\sigma)\) generated by the family of \(C\sigma(Ax+b)\) is dense in \(C(K, \mathbb{R})\). Since
\(V(\sigma)\subset V(ReLU)\),
it follows that \(V(ReLU)\) is also dense in \(C(K, \mathbb{R})\).

\end{proof}

We can deduce the following corollary:

\medskip

{\bf Corollary 5.3.}
Let \( (S, f) \) be a dataset, where \( S \subseteq \mathbb{R}^n \) consists of \( N \) elements, and \( f: S \to \underline{m} = \{1, 2, \dots, m\} \) is a function assigning each element in \( S \) to one of \( m \) classes.  There exists a neural network whose output is a function  \( \hat{Y}: \mathbb{R}^n \to \mathbb{R}^m \)  such that
\[
f(x) = \arg\max(\hat Y(x), \text{ axis}=1)
\]
for all \( x \in S \).

\medskip 

\begin{proof}
Let \( Y: S \to \mathbb{R}^m \) be the one-hot encoding of \( f \), i.e., for each \( x \in S \), \( Y(x) \) is the vector in \( \mathbb{R}^m \) with a 1 at position \( f(x) \) and 0 elsewhere.
The Whitney extension theorem implies that
 here exists a differentiable function \( Y': \mathbb{R}^n \to \mathbb{R}^m \) that extends  \( Y(x) \) over a compact neighbourhood of \( S \) in \( \mathbb{R}^n \), satisfying \( Y'(x) = Y(x) \) for all \( x \in S \).

 By the Universal Approximation Theorem, there exists a neural network \( \hat{Y}: \mathbb{R}^n \to \mathbb{R}^m \) such that for any \( \varepsilon > 0 \),
\[
\| \hat{Y}(x) - Y'(x) \| < \varepsilon \quad \text{for all } x \in S.
\]
 If \( \hat{Y}(x) \) is sufficiently close to \( Y'(x) \), the argmax of \( \hat{Y}(x) \) will be equal to the argmax of \( Y(x)\).

\end{proof}

While the existence of a neural network whose outputs can recover the classes of any labeled dataset follows from the universal approximation theorem, this result does not constitute a complete solution to the classification problem.

In particular, two major issues remain:

\bigskip

 {\bf Overfitting.}
 
\medskip
 
  A neural network that exactly recovers the labels on the training data may fail to generalize to unseen data. As shows the following example:

 Suppose that the \emph{rate of success} of a student is a linear function of the average number of hours he studies per week and his number of absences. Assume we are given $N$ students represented by data points $(x_1, y_1), (x_2, y_2), \ldots, (x_N, y_N)$, where each point corresponds to two features: study time and absences.

When plotted in the plane, the data is approximately linearly separable; that is, there exists a line which classifies most of the points correctly into ``successful'' and ``unsuccessful'' categories. However, due to noise or variability, a few outliers may be misclassified.

One may construct a high-degree polynomial curve that passes through all the data points and achieves perfect classification. Despite its low training error, such a model is likely to \emph{overfit} the training data and perform poorly on new, unseen examples.

\medskip

In order to avoid avoid overfitting, we can:

\medskip

{\bf Split the data.}

\medskip

We split the dataset sets into
two subsets: the training set and the validation set. The training set is used to
train the model: If the model performs significantly better on the training dataset than on the validation dataset, this suggests overfitting.

\medskip

{\bf Regularization.}

\medskip

A regularization term can be added to the loss function:

\(Lp\)-regularization:

A term proportional to the norm is added to the loss function:
\[
Loss_r(x)=Loss+\lambda\|x\|^p
\]

where \(\lambda \in \mathbb{R}\). This penalizes large weights.

\medskip

{\bf Dropout.}

Dropout is a regularization technique used in training neural networks. During each training epoch, a random percentage of neurons are temporarily deactivated (i.e., "dropped out"). This prevents the network from becoming overly reliant on specific neurons, thereby reducing overfitting. By forcing the network to operate without certain neurons during training, dropout encourages the development of independent and redundant feature representations, reducing interdependence (or co-adaptation) between neurons.

\bigskip

{\bf Complexity of decision boundaries.}

\medskip

Let \( S \) be a dataset embedded in the vector space \( \mathbb{R}^n \), and let \( f : S \rightarrow \{0, \dots, m-1\} \) be a labeling function. Suppose that \( S \) is clustered into \( m \) regions of \( \mathbb{R}^n \), separated by a finite collection of hypersurfaces. Since \( S \) is finite, we may assume these hypersurfaces are compact.

By a theorem of Nash (see [4] Theorem [3]), the decision boundaries separating these regions can be approximated arbitrarily closely by real algebraic varieties. Thus, we may suppose that the decision boundary is represented by an algebraic hypersurface \( H \subset \mathbb{R}^n \).

\medskip

Overfitting occurs when the decision boundary \( H \) is excessively complex (e.g., has high degree or intricate geometry). This motivates the following definition:

\medskip

\textbf{Definition 5.4.}  
Let \( r \in [0, 1] \) be a given accuracy threshold.  
A \emph{good decision boundary of order \( r \)} is an algebraic hypersurface \( H \subset \mathbb{R}^n \) that correctly classifies at least an \( r \)-fraction of the dataset \( S \), and among all such hypersurfaces, has the minimal possible algebraic degree.

\bigskip
  
 {\bf Computational feasibility:} 

 In practice, training neural networks   may be computationally expensive, for example the neural network used in the universal approximation theorem and the dataset maybe very large.
 This issue can be solved by:

 Splitting the dataset into batches of equal size. At each epoch, gradient descent is performed on each batch, and the loss for the epoch is calculated as the average of batch losses.
 
 Using graphical processing units (GPU) to accelerate computations.
 
 Creating specialized networks:
 
 For specific classification problems, more effective and computationally efficient linear maps can be used. For example, for image classification, convolutional layers perform very well.

\section*{6. Convolutional Neural Networks.}

Convolutional neural networks (CNN) have been created by Yann Le Cun in 1988 (see [10]). They where firstly applied to classify hand written digits. ImageNet is a dataset containing 1.2 millions images, 50000 validation images and 1000 classes.
The imageNet competition ran between 2010 and 2017, its goal was to find the algorithm which classify the ImageNet dataset with the best accuracy. CNNs gained widespread popularity when AlexNet won the 2012 ImageNet competition. AlexNet significantly outperformed all previous approaches, reducing the  error rate from 26 percent to 15.3 percent, and demonstrating the power of deep learning when combined with GPU acceleration.

A digital grayscale image can be represented as a 2D matrix, where each element corresponds to the intensity of a pixel. A color image is  represented by three 2D matrices, corresponding to the Red, Green, and Blue (RGB) color channels.

In neural network theory, a convolution is a  linear operation designed to exploit the spatial structure of images. It involves the use of kernels (also called filters), which are small matrices that slide over the input image. At each position, the kernel performs an element-wise multiplication followed by a summation, producing a single value in the output feature map. This operation allows the network to detect local patterns, such as edges, textures, or more complex features, depending on the layer depth.

More precisely:

\bigskip

\textbf{Definition 6.1.}

Let \( I \) be a \( C_{\text{in}} \times m \times n \) input tensor, where:
     \( C_{\text{in}} \) is the number of input channels,
    \( m \) and \( n \) are the spatial dimensions (height and width) of the input.

Let \( K \) be a \( C_{\text{out}} \times C_{\text{in}} \times p \times p \) filter (or kernel) tensor, where:
     \( C_{\text{out}} \) is the number of output channels,
    \( p \times p \) is the spatial size of the filter (height and width).
Let \( s \) be the \textbf{stride}.

The \textbf{convolution} of \( I \) by \( K \), denoted \( O \), is the output tensor with dimensions
\(
{C_{\text{out}} \times \left( \left\lfloor \frac{m - p}{s} \right\rfloor + 1 \right) \times \left( \left\lfloor \frac{n - p}{s} \right\rfloor + 1 \right)}
\)

and is defined by:

\[
O(d, i, j) = \sum_{c=0}^{C_{\text{in}} - 1} \sum_{u=0}^{p - 1} \sum_{v=0}^{p - 1} I(c, i s + u, j s + v) \cdot K(d, c, u, v)
\]

where
     \( i \) and \( j \) are the spatial indices in the output (height and width),
     \( d \) is the output channel index,
     \( c \) is the input channel index,
     \( u \) and \( v \) are the spatial indices in the filter (height and width),
     \( s \) is the stride.

\medskip

    The output size in terms of height and width is calculated based on the stride and filter size, and  \( \left\lfloor \cdot \right\rfloor \) is the floor operator defined on real numbers.
    
     The convolution operation sums over all input channels \( c \), all filter spatial positions \( u, v \), and computes the result for each spatial position \( i, j \) in the output tensor.
The stride \( s \) dictates how much the filter shifts for each convolution operation.

\medskip

{\bf Remarks.}

Convolutions maybe very cheaper than linear maps in the sense that they use less parameters:  suppose that \(C_{in}=C_{out}=1\), and \(I\) is an \(224\times 224\) image, usually a kernel is an \(p\times p\) matrix with \(p\leq 7\). which represents at most \(49\) parameters to which we add a parameter for the bias term.,  A classical linear layer with one neuron is defined by \(224\times 224\) parameters and \(1\) parameters for the biase.

\medskip

To reduce memory usage, convolutional neural networks are often trained in batches by dividing the dataset into smaller, equally sized parts that fit into the available RAM.

The following layers are often added to convolutional layers:

\medskip

{\bf Padding.}

\medskip

Padding may be added to the original image to increase the dimensions of the output feature map after convolution.  
More precisely, a \emph{padding of size} \( d \) applied to an \(C_{in}\times  m \times n \) matrix \( I \) consists of embedding \( I \) into a new matrix \(C_{in}\times (m+2d)\times (n+2d)\) matrix \(J\) defined by: 
\[
J(c,i, j) =
\begin{cases}
I(c,i - d, j - d) & \text{if } d \leq i < m + d \text{ and } d \leq j < n + d, \\
0 & \text{otherwise}.
\end{cases}
\]

{\bf Pooling.}

\medskip

Pooling can also be added to convolutional neural networks to reduce the size of the output image. It is defined as follows:

Max Pooling

Let \(I\) be an \(C_{in}\times m\times n\) matrix, \(p\) the size of the pooling window, and \(s\) the stride. The output of the pulling of \(I\) with the kernel \(K\) is the 
\( (\left\lfloor \frac{m - p}{s} \right\rfloor + 1)\times (\left\lfloor \frac{n - p}{s} \right\rfloor + 1)\)
matrix defined by:
\[
O(c,i,j)=\max_{0\leq u<p,0\leq v<p}(c,is+u,js+v)
\]

Average Pooling

For average pooling, the output image is defined by:
\[
O(c,i,j)=\frac{1}{p^2}\sum_{u=0}^{p-1}\sum_{v=0}^{p-1}I(c,si+u,sj+v)
\]

\medskip

{\bf Batch Normalization.}

It is  a technique used in neural networks to normalize the inputs of each layer, stabilizing and speeding up training by reducing internal covariate shift.

It is defined as follows: Let \(x_i\) be the \(i\)-sample of the batch. We denote by \(\mu_B\) the mean of the batch, and by \(\sigma_B\) its variance, we have:

\begin{align*}
\hat{x}_i &= \frac{x_i - \mu_B}{\sqrt{\sigma_B^2 + \epsilon}} \\
y_i &= \gamma \hat{x}_i + \beta
\end{align*}

where 
 \( \hat{x}_i \) is the  normalized input
   \( \epsilon \), a small constant for numerical stability (avoiding division by zero)
  \( \gamma \), a learnable (updapted through gradient descent) scaling parameter
   \( \beta \), learnable shifting parameter
  \( y_i \): final output after normalization.

\medskip

{\bf Dropout}

It is  a regularization method where randomly selected neurons are ignored during training, helping prevent overfitting by making the network less reliant on specific neurons.

\subsection*{Training a convolution neural network.}

The training process for CNNs follows a scheme similar to that used for fully connected deep neural networks. In particular, training is done via backpropagation.
The forward pass of a CNN typically starts with several convolutional blocks. Each convolutional block usually consists of:

\begin{itemize}
    \item A \textbf{convolutional layer}, which applies learnable filters (kernels) to the input. These are \textit{linear maps}.

    \item A \textbf{normalization or regularization layer}, such as \textit{Dropout} or \textit{Batch Normalization}, to prevent overfitting and stabilize the learning process.
    \item A \textbf{pooling layer} (e.g., MaxPooling or AveragePooling) to reduce spatial dimensions and computation.
\end{itemize}

After several convolutional blocks, the output tensor is passed through a Flatten layer that converts the multi-dimensional output into a one-dimensional vector. This vector is then fed into one or more fully connected (linear) layers, often used for classification or regression tasks.

Note that \textbf{convolutional layers are linear operators}. Therefore, convolutional neural networks can be trained using gradient-based optimization techniques (e.g., Stochastic Gradient Descent, Adam), with backpropagation used to compute the gradients.

\subsection*{Using Software libraries.}

In the early stages of deep learning research and development, neural network models were often implemented entirely from scratch. In the case of Python, developers commonly used \texttt{NumPy} to handle tensor operations. This required manual implementation of the forward pass, backward pass (backpropagation), optimization routines (such as gradient descent), and loss functions. While this approach helped researchers gain a deep understanding of the underlying mechanics, it was time-consuming and error-prone, especially for large and complex models.

The emergence of automatic differentiation libraries significantly streamlined this process. The first major breakthrough was \textbf{Theano}, which introduced symbolic computation and automatic differentiation, enabling users to define models and compute gradients without manually deriving formulas.

Building on this foundation, newer and more user-friendly libraries such as \textbf{TensorFlow} and \textbf{PyTorch} further popularized deep learning by providing high-level APIs, prebuilt layers (e.g., convolutional, linear, pooling), loss functions, and optimizers. These libraries automatically handle backpropagation through computational graphs and offer tools for efficient GPU computation, dataset management, and model deployment.

Today, using such APIs has become the standard practice in both research and industry. They allow developers to focus on designing and experimenting with model architectures rather than implementing low-level numerical routines, thus accelerating the development cycle and enabling rapid prototyping of deep learning solutions.

\medskip

{\bf The malaria dataset.}

\medskip

We have constructed a convolutional neural network using PyTorch to classify blood cells infected by malaria. For this purpose, we used the malaria dataset available in the TensorFlow library.

We split the dataset into \(80\%\) for training and \(20\%\) for validation.

We designed a convolutional neural network composed of three convolutional blocks followed by four fully connected layers.

The first convolutional block consists of:
a convolutional layer with \(16\) filters and a kernel size of \(3 \times 3\),
a ReLU activation function,
and a max pooling layer with a \(2 \times 2\) window.

The second convolutional block consists of:
a convolutional layer with \(32\) filters and a kernel size of \(3 \times 3\),
a ReLU activation function,
and a max pooling layer with a \(2 \times 2\) window.

The third convolutional block consists of:
a convolutional layer with \(64\) filters and a kernel size of \(3 \times 3\),
a ReLU activation function,
and a max pooling layer with a \(2 \times 2\) window.

These convolutional blocks are followed by four linear layers with the following sizes:
\(4096 \times 512\),
\(512 \times 256\),
\(256 \times 128\),
and \(128 \times 2\).
The first three linear layers use ReLU activation.

We used the Adam optimizer and the cross-entropy loss function.

The network was trained for 5 epochs and achieved a validation accuracy of over \(95\%\). This high accuracy is likely due to the relatively simple structural features of blood cells.

We have deployed this model in the cloud with gradio.

\bigskip

{\bf Pneumonia dataset.}

\medskip

The malaria dataset consists of simple elements, allowing us to build a model that achieves good classification accuracy. However, when datasets contain more complex features, this becomes more challenging. In such cases, fine-tuning is commonly used: a pretrained model from a related task is adapted by replacing its final layers to suit the current task. This approach, known as transfer learning, is supported by libraries like HuggingFace,Tensorflow, Pytorch,  which provide many freely available pretrained models.

Pneumonia kills around {2 million children under 5 years old} globally each year. One common method of diagnosis is through chest X-rays. However, in many regions, there is a shortage of specialists capable of interpreting radiologic images. {Deep learning} offers a potential solution.

Kermany \textit{et al. see [9]} applied {transfer learning} to develop a pneumonia detection algorithm, fine-tuning the \textbf{Inception V3} model and achieving an accuracy of \textbf{92.8\%}.

In our study, we used the same dataset but fine-tuned {EfficientNet B0}, achieving an improved accuracy of \textbf{93.5\%} after just {5 training epochs}. To achieve this, we applied a data augmentation technique called {Mixup}, which creates new training examples by interpolating both images and their labels. The model was trained on this augmented data and evaluated on the original, unseen validation set, resulting in a final accuracy of \textbf{93.5\%}.

 The code is available in the Github repository.

\bigskip

\section*{7. Recurrent Neural Network.}

We have seen that convolutional neural networks (CNNs) are a class of neural networks specifically designed to process data with a grid-like topology, such as images. 
 Recurrent neural networks (RNNs) are neural networks designed to handle {sequential data}, these are data where the order of elements is relevant such as time series, natural language text, or audio signals. These models are structured to recognize patterns and dependencies over time, making them well-suited for tasks involving sequences of varying lengths.

\medskip

\noindent
\textbf{Simple RNN}

\medskip

To illustrate the idea behind recurrent neural networks, let us consider the architecture of a simple RNN. A simple RNN consists of three main components:

\begin{itemize}
    \item \textbf{Input neurons} (which process elements of the input sequence),
    \item \textbf{Hidden neurons} (which maintain a memory of previous inputs through recurrent connections),
    \item \textbf{Output neurons} (which produce predictions  based on the current hidden state).
\end{itemize}

The defining feature of an RNN is the presence of a \textit{loop} (recurrent connection) in the hidden layer: each hidden neuron receives not only the current input but also its own output from the previous time step. This loop allows the network to maintain a form of memory, capturing information from earlier in the sequence and using it to influence later outputs.

\medskip

Represent an input sequence by \(x_0,\cdots x_{T-1}\).
Firstly, \(x_0\) is fed in the unit, we suppose that the initial hidden state \(h_{-1}\) is \(0\). The second hidden state is:

\[
h_1 = \sigma(W_{xh} x_0 + W_{hh} h_{-1} + b_h)
\]

and the first output is:
\[
y_1 = \phi(W_{hy} h_1 + b_y)
\]

where
     \( W_{xh}, W_{hh}, W_{hy} \) are weight matrices,
    \( b_h, b_y \) are biases,
     \(\sigma,\psi\): are activation functions.

Recursively, we define \(h_t\) and \(x_t\).

\[
h_t = \sigma(W_{xh} x_t + W_{hh} h_{t-1} + b_h)
\]
\[
y_t = \phi(W_{hy} h_t + b_y)
\]

The network produces a series of outputs \(y_1,\cdots y_T\) such that \(y_t\) depends on \(y_{t-1}\). Usually, the activation of the hidden node is the hyperbolic tangent.

\medskip

    
    
    

{\bf  Gradient of loss with respect to hidden states.}

\medskip

Assume the total loss is:
\[
\mathcal{L} = \sum_{t=1}^{T} \mathcal{L}_t(h_t)
\]

We want to compute the gradient \( \frac{\partial \mathcal{L}}{\partial h_t} \), which depends on all future steps \( k \geq t \). By the chain rule:

\[
\frac{\partial \mathcal{L}}{\partial h_t} = \sum_{k=t}^{T} \frac{\partial \mathcal{L}_k}{\partial h_k} \cdot \frac{\partial h_k}{\partial h_t}
\]

where the output sequence is \(\hat y_0(h_0),\cdots,\hat y_{T-1}(h_{T-1})\), the expected output is \(y_0,...,y_{T-1}\) and the loss at each time step is  \(\mathcal{L}_t(h_t)=L(\hat y_t(h_t),y_t)\).
\medskip

To compute \( \frac{\partial h_k}{\partial h_t} \), consider:

\[
\frac{\partial h_j}{\partial h_{j-1}} = \frac{\partial h_j}{\partial a_j} \cdot \frac{\partial a_j}{\partial h_{j-1}} = \text{diag}(\phi'(a_j)) \cdot W_{hh}
\]
where:
\[
a_j = W_{xh} x_j + W_{hh} h_{j-1} + b_h
\]

and \(diag(\phi'(a_j))\) is a diagonal matrix containing the derivative of \(\phi\) with respect to \(a_j\) applied element wise.

Thus, for \( k > t \):

\[
\frac{\partial h_k}{\partial h_t} = \prod_{j=t+1}^{k} \left( \text{diag}(\phi'(a_j)) \cdot W_{hh} \right)
\]

\medskip

The full gradient with respect to the hidden state at time \( t \) is:

\[
\frac{\partial \mathcal{L}}{\partial h_t} = \sum_{k=t}^{T} \left( \frac{\partial \mathcal{L}_k}{\partial h_k} \cdot \prod_{j=t+1}^{k} \text{diag}(\phi'(a_j)) \cdot W_{hh} \right)
\]

This recursive structure is the source of the vanishing and exploding gradient problem, since repeated multiplication of Jacobians can shrink or grow exponentially depending on the spectral norm of \( W_{hh} \) and the activation derivative \( \phi'(a_j) \);
 LSTM or GRU can solve this issue.

\bigskip

{\bf Long Short-Term Memory (LSTM)}

\bigskip

Long Short-Term Memory (LSTM) networks are a special kind of recurrent neural network  capable of learning long-term dependencies. They introduce a cell state and a gate state to control information flow.

\medskip

{\bf LSTM Equations}

\medskip

Given input \( x_t \), previous hidden state \( h_{t-1} \), and previous cell state \( c_{t-1} \):

\begin{align*}
f_t &= \sigma(W_f x_t + U_f h_{t-1} + b_f) && \text{(Forget gate)} \\
i_t &= \sigma(W_i x_t + U_i h_{t-1} + b_i) && \text{(Input gate)} \\
\tilde{c}_t &= \tanh(W_c x_t + U_c h_{t-1} + b_c) && \text{(Candidate cell state)} \\
c_t &= f_t * c_{t-1} + i_t * \tilde{c}_t && \text{(New cell state)} \\
o_t &= \sigma(W_o x_t + U_o h_{t-1} + b_o) && \text{(Output gate)} \\
h_t &= o_t * \tanh(c_t) && \text{(New hidden state)}
\end{align*}

Where:
\( \sigma \) is the sigmoid function,
 \( \tanh \) is the hyperbolic tangent function, and
 \( * \) denotes the Hadamard product of matrices, that is the element-wise multiplication.

The input gate controls how much new information from the current input should be added to the memory.
The forget gate controls how much of the past cell state \(c_{t-1}\) should be retained or forgotten.
The cell state carry the long memory over states.
It seems that LSTM can preserve gradients across long time intervals, which mitigates the vanishing gradient problem in standard RNNs. We have not analyzed the expression of backpropagation from a mathematical perspective to show that they converge faster than simple RNN.

\bigskip

{\bf Gated Recurrent Unit (GRU)}

\bigskip

GRU is a simplified version of LSTM that combines the forget and input gates into a single update gate, and merges the cell and hidden states.

\medskip

{\bf GRU Equations}

\medskip

Given input \( x_t \), and previous hidden state \( h_{t-1} \):

\begin{align*}
z_t &= \sigma(W_z x_t + U_z h_{t-1} + b_z) && \text{(Update gate)} \\
r_t &= \sigma(W_r x_t + U_r h_{t-1} + b_r) && \text{(Reset gate)} \\
\tilde{h}_t &= \tanh(W_h x_t + U_h (r_t * h_{t-1}) + b_h) && \text{(Candidate activation)} \\
h_t &= (1 - z_t) * h_{t-1} + z_t * \tilde{h}_t && \text{(New hidden state)}
\end{align*}

It seems that GRUs performs similarly to LSTMs but are computationally more efficient due to fewer parameters.

\subsection*{Attention Mechanism.}

In natural language processing tasks such as machine translation, the meaning of a word in a sentence is influenced by the meanings of the surrounding words. This is particularly important when dealing with homonyms.

Consider the following two sentences:
\begin{itemize}
    \item \textit{He went to the bank to get money.}
    \item \textit{He sat by the bank of the river.}
\end{itemize}

The word \textit{bank} has different meanings in each sentence. The attention mechanism introduced in the transformer model provides a mathematical framework to capture such contextual relationships.

\medskip

{\bf Mathematical Formalization}

\medskip

Assume that each word is represented by an embedding vector \( x \in \mathbb{R}^d \). For each word, we compute:

\begin{align*}
    x^q &= W^Q x \in \mathbb{R}^{d_k} \quad \text{(query vector)} \\
    x^k &= W^K x \in \mathbb{R}^{d_k} \quad \text{(key vector)} \\
    x^v &= W^V x \in \mathbb{R}^{d_l} \quad \text{(value vector)}
\end{align*}

where \( W^Q, W^K \in M_{d\times d_k}(\mathbb{R}) \) and \( W^V \in M_{d\times d_l}(\mathbb{R}) \) are learnable weight matrices.

For a pair of words \( x_i \) and \( x_j \), the  attention score is given by:

\[
a_{ij} = \frac{\langle W^Q x_i, W^K x_j \rangle}{\sqrt{d_k}} = \frac{\langle x_i^q, x_j^k \rangle}{\sqrt{d_k}}
\]

The attention weights are:

\[
\alpha_{i1}, \ldots, \alpha_{ip} = \text{softmax}(a_{i1}, \ldots, a_{ip})
\]

The output vector (contextual representation) of word \( x_i \) is:

\[
z_i = \sum_{j=1}^{p} \alpha_{ij} W^V x_j = \sum_{j=1}^{p} \alpha_{ij} x_j^v
\]

This representation \( z_i \) depends on the entire sentence and varies with context, even for the same input word vector \( x_i \).

\section{Transformer}

The transformer architecture  which is the core of large language models today,   processes sequences using attention mechanisms instead of recurrence. We present a version of the transformer with one attention head.

Let the input sequence be \( X = [x_1, x_2, \dots, x_n] \), where each \( x_i \in \mathbb{R}^d \) is a word embedding. We can define the Query, Key, and Value matrices with the following formulas:

\begin{itemize}
    \item \textbf{Query: } \( Q = X W^Q \in M_{n\times d_k}(\mathbb{R}) \)
    \item \textbf{Key: } \( K = X W^K \in M_{n\times d_k}(\mathbb{R}) \)
    \item \textbf{Value: } \( V = X W^V \in M_{n\times d_v}(\mathbb{R}) \)
\end{itemize}

Where:

\begin{itemize}
    \item \( W^Q, W^K \in M_{d\times d_k}(\mathbb{R}) \)
    \item \( W^V \in M_{d\times d_v}(\mathbb{R}) \)
\end{itemize}

We compute the attention scores and outputs as follows:

\[
A = \text{softmax}\left( \frac{Q K^\top}{\sqrt{d_k}} \right) \in M_{n\times n}(\mathbb{R})
\]

\[
Z = A V \in M_{n\times d_v}(\mathbb{R})
\]

Here, \( Z \) is a context-aware representation of the input sequence.

We then apply a linear layer to each row \( z_i \) of \( Z \):

\[
\text{FFN}(z_i) = \text{ReLU}( z_i W_1 + b_1) W_2 + b_2
\]

\medskip

\textbf{Residual and Normalization Layer.}

\medskip

After applying the feed-forward network (FFN) to each row, the output is added to the original input in a residual connection. This can be expressed as:

\[
\text{Res} = X + \text{FFN}(Z)
\]

Finally, a layer normalization is applied to the result of the residual connection:

\[
\text{Output} = \text{LayerNorm}(\text{Res})
\]

This ensures stable training and improves gradient flow.

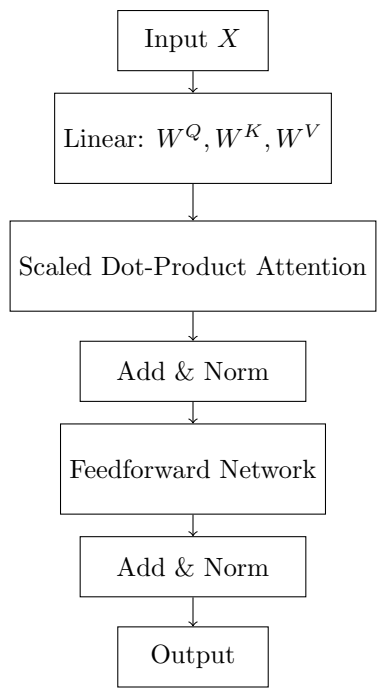
\begin{center}
\begin{tikzpicture}[scale=1, every node/.style={scale=1}]
\node (x) at (0, 0) [draw, minimum width=2cm, minimum height=0.8cm] {Input $X$};

\node (qkv) at (0, -1.3) [draw, minimum width=3cm, minimum height=1.2cm] {Linear: $W^Q, W^K, W^V$};

\node (attn) at (0, -3) [draw, minimum width=3.5cm, minimum height=1.2cm] {Scaled Dot-Product Attention};

\node (an1) at (0, -4.4) [draw, minimum width=3cm, minimum height=0.8cm] {Add \& Norm};

\node (ffn) at (0, -5.7) [draw, minimum width=3cm, minimum height=1.2cm] {Feedforward Network};

\node (an2) at (0, -7) [draw, minimum width=3cm, minimum height=0.8cm] {Add \& Norm};

\node (out) at (0, -8.2) [draw, minimum width=2cm, minimum height=0.8cm] {Output};

\draw[->] (x) -- (qkv);
\draw[->] (qkv) -- (attn);
\draw[->] (attn) -- (an1);
\draw[->] (an1) -- (ffn);
\draw[->] (ffn) -- (an2);
\draw[->] (an2) -- (out);
\end{tikzpicture}
\end{center}

With this structures, the neural network learns to focus on important context tokens. Each word's output representation is computed based on its relationships (via dot products) with all other words, giving rise to rich contextual embeddings.

\bigskip

{\bf Training Neural network for Natural language processing tasks.}

\bigskip

In a natural language processing task, the dataset typically consists of raw text. Before training, the text is processed by a tokenizer that splits words into smaller units called tokens, often using algorithms such as Byte Pair Encoding (BPE), WordPiece, or other subword segmentation methods. These tokens are then mapped into a high-dimensional vector space through an embedding layer, which is trained alongside the model to capture semantic and syntactic relationships between words. The embedded token sequences are input into a neural network—such as a Transformer—which learns to perform the desired task. Training is guided by a suitable loss function, such as cross-entropy loss for language modeling or classification, enabling the model to optimize its performance through backpropagation. In large language models (LLMs), alignment is performed after pretraining to ensure that the model’s behavior aligns with human values, intentions, and preferences. This process typically involves fine-tuning the model using techniques such as supervised learning from human-annotated examples and reinforcement learning from human feedback (RLHF).

It is clear that the attention mechanism has shown its effictiveness. It is at the heart of large language models, like GPT, LLAMA, DeepSeek,... We have not analyze from a mathematical perspective why it converges faster than RNN. 

\section*{8. The category of neural networks.}

{\bf Graph neural networks and homotopy.}

\medskip

In this section, we study \textbf{Graph Neural Networks (GNNs)} — a general framework that unifies deep neural networks and recurrent neural networks over graph-structured data.

\medskip

\textbf{Definition 8.1.} A \emph{directed graph} \( G \) is defined by a set of nodes \( N_G \), a set of arcs \( A_G \), and source and target maps:
\[
s, t : A_G \rightarrow N_G
\]

It can also be viewed as a \emph{topos of sheaves} on the category \( C_G \), which has two objects \( 0, 1 \) and two morphisms:
\[
s_G, t_G : 0 \rightarrow 1
\]

This alternative view gives rise to a category \( \mathbf{Gph} \), where morphisms are morphisms between presheaves. Precisely, given graphs \( G = (N_G, A_G) \) and \( G' = (N_{G'}, A_{G'}) \), a morphism \( f: G \rightarrow G' \) consists of two maps:
\[
f_0: N_G \rightarrow N_{G'}, \quad f_1: A_G \rightarrow A_{G'}
\]
such that:
\[
f_0 \circ s = s' \circ f_1 \quad \text{and} \quad f_0 \circ t = t' \circ f_1
\]

Let \(x\) be a node of the graph \(G\), we denote by \(G(*,x)\) the set of arrows of \(G\) whose target is \(x\), and by \(G(x,*)\) the set of arcs of \(G\) whose source is \(x\).

\bigskip

{\bf Definitions 8.2.}

\medskip

Let \( n \geq 1 \) be an integer. The \emph{\(n\)-cycle graph} \( C_n \) is the graph with node set:
\[
N_{C_n} = \{0, 1, \ldots, n-1\}
\]
There exists a unique arc from \( i \) to \( i+1 \) for each \( i < n-1 \), and a unique arc from \( n-1 \) to \( 0 \), forming a closed directed cycle.

A graph \(G\) is acyclic if for every integer \(n\geq 1\), the set of morphisms \(Gph(C_n,G)\) is empty. 

\medskip

\textbf{Definition 8.3.}

A \emph{Graph Neural Network (GNN)} is defined by the following data:

\begin{itemize}
    \item A directed graph \( G \), together with a collection of subsets of its nodes:
    \[
    L = \{N_0, N_1, \ldots, N_p\}
    \]

    \item For every node \( x_{i+1} \in N_{i+1} \), there exists a node \( x_i \in N_i \) and an  arc:
    \[
    a_{x_i, x_{i+1}} : x_i \rightarrow x_{i+1}
    \]
    in \( G \).

    \item A vector space \( E_x \) is associated to each node \( x \in G \).
\end{itemize}

\medskip

\textbf{Information circulating through a GNN.}

\medskip

The GNN updates node information over discrete time steps. Suppose the initial information is given by a family of elements \( (u_x)_{x \in N_0} \), where \( u_x \in E_x \).

At timestep \( t \), the current layer is \( N_{r_t} \), where:
\[
r_t = t \bmod (p + 1)
\]
For every node \( x \in N_{r_t} \), the feature vector at time \( t \) is \( u_x \in E_x \).

At timestep \( t+1 \), we update the features of nodes in the layer \( N_{r_{t+1}} \), where \( r_{t+1} = (t+1) \bmod (p+1) \).

Let \( x \in N_{r_{t+1}} \). Define \( A^N_x \) as the set of arcs \( a: y \rightarrow x \) in \( G \), where \( y \in N_{r_t} \). For each such \( y \), we assume there is a linear map:
\[
l_{yx} : E_y \rightarrow E_x
\]

The updated information at node \( x \) is then given by:
\[
u_x = \sigma\left( \sum_{a \in A^N_x} l_{yx}(u_y) \right)
\]
where \( \sigma \) is a fixed non-linear activation function.
  
\bigskip

{\bf Remarks.}

\medskip

In the classical definition of Graph Neural Networks (GNNs), it is assumed that the graph has only one layer.

Deep neural networks can be viewed as a special case of GNNs where the underlying graph is acyclic.

A simple recurrent network corresponds to a GNN defined by a graph with node set \( G_0 = \{0, 1, 2\} \). There exists a unique arc from node \( 0 \) to node \( 1 \), a self-loop at node \( 1 \), and a unique arc from node \( 1 \) to node \( 2 \).

\medskip

The fact that the computational graph associated with a simple RNN is {not acyclic}, while the graph associated with a DNN {is acyclic}, suggests that {cycles in computational graphs are related to the ability of the corresponding neural network to process memory}.

In previous work, we defined a closed model structure in the topos of graphs, where the class of weak equivalences $\mathcal{W}$ consists of graph morphisms that induce {bijections between cycles}. Another model, with the same class of weak fibrations, was defined in the work of Tsemo Bisson see( [1] and [13]).

This suggests that the {homotopy class of a graph}, as defined by this model, can be interpreted as a representation of its {memory capacity}.

\bigskip

{\bf Tensor categories and graph neural networks.}

\medskip

Let \(G\) a directed graph, for every node \( x \) of \( G \), 
the set of initial nodes \( G_{\text{in}} \) of \( G \) is the subset of nodes of \( G \) such that for every \( x \in G_{\text{in}} \), \( G(*, x) = \emptyset \).

The set of final nodes \( G_{\text{out}} \) of \( G \) is the subset of nodes of \( G \) such that for every \( x \in G_{\text{out}} \), \( G(x, *) = \emptyset \).

We now define the category \( \mathbf{Net} \) whose objects are finite sets. A morphism \( f: A \to B \) between objects of \( \mathbf{Net} \) is an oriented graph \( F \) such that \( F_{\text{in}} = A \) and \( F_{\text{out}} = B \).

Let \( f: A \to B \) and \( g: B \to C \) be two morphisms of \( \mathbf{Net} \), represented by the oriented graphs \( F \) and \( G \), respectively. The composition \( g \circ f: A \to C \) is represented by the graph whose set of nodes is \( F_0 \cup G_0 \) and whose set of arcs is \( F_1 \cup G_1 \). We denote by \( \text{Hom}_{\mathbf{Net}}(A, B) \) the set of morphisms from \( A \) to \( B \).

The identity morphism in \( \text{Hom}_{\mathbf{Net}}(A, A) \) is the graph whose set of nodes is \( A \) and whose set of arcs is empty.

\bigskip

The category \( \mathbf{Net} \) is equipped with a tensor structure, defined as follows:

- For objects \( A \) and \( B \) in \( \mathbf{Net} \), the tensor product \( A \otimes B \) is the disjoint union of \( A \) and \( B \).

- For every morphisms \( f: A \to A' \) and \( g: B \to B' \) in \( \mathbf{Net} \), the tensor product \( f \otimes g: A \otimes B \to A' \otimes B' \) is the disjoint union of the graphs \( f \) and \( g \).

\bigskip
\bigskip
\bigskip

\section*{Bibliography.}

\bigskip
\bigskip

1. Bisson, T.,  Tsemo, A. (2009). A homotopical algebra of graphs related to zeta series.

\medskip

2. Bahdanau, D., Cho, K., and Bengio, Y. (2015). Neural machine translation by jointly learning to align and translate. Proceedings of the International Conference on Learning Representations (ICLR).

\medskip
3. Chollet, F. (2021). Deep learning with Python. simon and schuster.
\medskip

4. Cybenko, G. (1989). Approximation by superpositions of a sigmoidal function. Mathematics of control, signals and systems, 2(4), 303-314.

\medskip

5. Fairchild M, Lemoine M. Topological data analysis of mortality patterns during the COVID-19 pandemic. https://arxiv.org/pdf/2510.15066

\medskip

6. Gatys, L. A., Ecker, A. S.,  Bethge, M. (2015). A neural algorithm of artistic style. arXiv preprint arXiv:1508.06576.

\medskip

7. Gowers, T - Why are LLMs not Better at Finding Proofs ?

$https://www.youtube.com/watch?v=5D3x_Ygv3No$

\medskip

8. Kassel, C. (2012). Quantum groups (Vol. 155). Springer Science  Business Media.

\medskip

9. Kermany, D. S., Goldbaum, M., Cai, W., Valentim, C. C., Liang, H., Baxter, S. L., ...  Zhang, K. (2018). Identifying medical diagnoses and treatable diseases by image-based deep learning. cell, 172(5), 1122-1131.

\medskip

10. LeCun, Y., Touresky, D., Hinton, G.,  Sejnowski, T. (1988, June). A theoretical framework for back-propagation. In Proceedings of the 1988 connectionist models summer school (Vol. 1, pp. 21-28).

\medskip

11. McCulloch, W. S.,  Pitts, W.  A logical calculus of the ideas immanent in nervous activity. The Bulletin of Mathematical Biophysics, 5, 115–133. (1943)

\medskip

12. Tognoli, A. (2006). Algebraic approximation of manifolds and spaces. In Séminaire Bourbaki vol. 1979/80 Exposés 543–560: Avec table par noms d'auteurs de 1967/68 à 1979/80 (pp. 73-94

\medskip

13. Tsemo A. (2017), Applications of closed models defined by counting to graph theory and topology. Algebra Letters 
\medskip

14. Vaswani, A., Shazeer, N., Parmar, N., Uszkoreit, J., Jones, L., Gomez, A. A., Kaiser, Ł., and Polosukhin, I. Attention is all you need. Proceedings of the Neural Information Processing Systems (NeurIPS), 30. (2017)

\medskip

15. Viola, P.,  Jones, M. (2001). Rapid object detection using a boosted cascade of simple features. In Proceedings of the 2001 IEEE computer society conference on computer vision and pattern recognition. CVPR 2001 (Vol. 1, pp. I-I). Ieee.

\medskip

16. Yang Huie  He. Mathematics, the rise of machines.

https://www.youtube.com/watch?v=oOYcPkBaotg

\end{document}